\documentclass[11pt]{article}
\usepackage{graphicx}
\usepackage{float}
\usepackage{subfig}
\usepackage{latexsym,amsmath,amsfonts,amscd, amsthm, dsfont}
\usepackage{nameref}

\usepackage{algorithm}
\usepackage{algpseudocode}
\usepackage{amsfonts}

\usepackage[percent]{overpic}
\usepackage{cleveref}
\usepackage{bm,color}
\usepackage{epsfig,verbatim,epstopdf,graphics}
\usepackage{changebar}
\usepackage{multirow}
\usepackage{enumerate}
\usepackage{yhmath}
 \usepackage{booktabs} 
 \usepackage{tikz}
\usepackage{verbatim}
\usepackage{appendix}
\usetikzlibrary{arrows,backgrounds,snakes,shapes}

\usepackage{soul}

\usepackage[noblocks]{authblk}

 \numberwithin{equation}{section}

\graphicspath{{./}{./figure/}}
\allowdisplaybreaks

\newtheoremstyle{plainNoItalics}{}{}{\normalfont}{}{\bfseries}{.}{ }{}

\theoremstyle{plain}
\newtheorem{thm}{Theorem}[section]
\theoremstyle{plainNoItalics}
\newtheorem{cor}[thm]{Corollary}

\newtheorem{prop}[thm]{Proposition}
\newtheorem{exa}[thm]{Example}

\DeclareMathOperator*{\argmax}{arg\,max}

\newcommand{\be}{\begin{eqnarray}}
\newcommand{\ee}{\end{eqnarray}}
\newcommand{\beno}{\begin{eqnarray*}}
\newcommand{\eeno}{\end{eqnarray*}}

\makeatletter

\newcommand{\Rmnum}[1]{\expandafter\@slowromancap\romannumeral #1@}
\makeatother

\title{A Semi-Lagrangian Adaptive-Rank (SLAR) Method for Linear Advection and Nonlinear Vlasov-Poisson System}

\author[a]{Nanyi Zheng}
\author[a]{Daniel Hayes}
\author[b]{Andrew Christlieb}
\author[a]{Jing-Mei Qiu}

\affil[a]{\small{Department of Mathematical Sciences, University of Delaware, Newark, DE 19716}}
\affil[b]{\small{Computational Mathematics, Science and Engineering, Michigan State University, East Lansing, 48824}}

\begin{document}
\maketitle







\bigskip
\noindent
{\bf Abstract.}
High-order semi-Lagrangian methods for kinetic equations have been under rapid development in the past few decades. In this work, we propose a semi-Lagrangian adaptive rank (SLAR) integrator in the finite difference framework for linear advection and nonlinear Vlasov-Poisson systems without dimensional splitting. The proposed method leverages the semi-Lagrangian approach to allow for significantly larger time steps while also exploiting the low-rank structure of the solution. This is achieved through cross approximation of matrices, also referred to as CUR or pseudo-skeleton approximation, where representative columns and rows are selected using specific strategies. To maintain numerical stability and ensure local mass conservation, we apply singular value truncation and a mass-conservative projection following the cross approximation of the updated solution. The computational complexity of our method scales linearly with the mesh size $N$ per dimension, compared to the $\mathcal{O}(N^2)$ complexity of traditional full-rank methods per time step. The algorithm is extended to handle nonlinear Vlasov-Poisson systems using a Runge-Kutta exponential integrator. Moreover, we evolve the macroscopic conservation laws for charge densities implicitly, enabling the use of large time steps that align with the semi-Lagrangian solver. We also perform a mass-conservative correction to ensure that the adaptive rank solution preserves macroscopic charge density conservation. To validate the efficiency and effectiveness of our method, we conduct a series of benchmark tests on both linear advection and nonlinear Vlasov-Poisson systems. The propose algorithm will have the potential in overcoming the curse of dimensionality for beyond 2D high dimensional problems, which is the subject of our future work. 




\section{Introduction}

High-order semi-Lagrangian (SL) methods have been developed over the past few decades to address both linear advection equations and nonlinear kinetic systems. The SL approach offers a unique blend of advantages from both Eulerian and Lagrangian perspectives. By using characteristic tracing, similar to the pure Lagrangian method, SL methods enable large time stepping. Simultaneously, they rely on a fixed Eulerian mesh, facilitating high-order spatial accuracy. Depending on the application, these methods can be designed using various spatial discretization techniques,
including finite element methods \cite{russell2002overview, carpio2014anisotropic, puigferrat2021semi}, discontinuous Galerkin (DG) methods \cite{restelli2006semi, rossmanith2011positivity, qiu2011positivity, cai2018high}, finite difference (FD) methods \cite{carrillo2007nonoscillatory, cristiani2007fast, xiong2019conservative}, and finite volume (FV) methods \cite{phillips2001semi, filbet2001conservative, lauritzen2010conservative}. Many existing SL methods use operator-splitting techniques due to their relative simplicity in handling high-dimensional problems \cite{carrillo2007nonoscillatory, qiu2011positivity}. However, many of such methods introduce splitting errors that can become significant in nonlinear cases. In contrast, non-splitting SL methods avoid these errors but present greater challenges in design and implementation, particularly when aiming to ensure local mass conservation. 

Parallel to the advancements in SL methods, the low-rank tensor approach has gained prominence as an effective strategy for addressing the curse of dimensionality and expediting high-dimensional kinetic simulations in recent decades. Such approach achieves compression of the numerical solution through low-rank decompositions of matrices and tensors. Two primary strategies have been developed for evolving low-rank solutions in time-dependent problems: the dynamic low-rank (DLR) approach \cite{koch2007dynamical, ceruti2022unconventional, einkemmer2018low} and the step-and-truncate (SAT) approach, also known as 'adaptive rank' methods \cite{kormann2015semi, dektor2021rank, GuoVlasovFlowMap2022, guo2024local, sands2024high}. For both linear transport and nonlinear kinetic models, the dynamic low-rank approach evolves solutions on a low-rank manifold using tangent space projections \cite{einkemmer2018low, einkemmer2020low}. In contrast, the step-and-truncate approach adapts the rank of the solution dynamically through tensor decomposition, manipulation, and truncation procedures, from direct discretization of partial differential operators \cite{oseledets2011tensor, kressner2014algorithm}. Under the step-and-truncate framework, significant advancements have been made, including semi-Lagrangian methods based on dimensional splitting \cite{kormann2015semi}, Eulerian methods that evolve solutions via a traditional method-of-lines approach in a low-rank format \cite{dektor2021rank, GuoVlasovFlowMap2022, sands2024high}, and collocation type low rank methods that leverage effective matrix sampling strategies  \cite{ghahremani2024cross, dektor2024collocation}. However, to date, no adaptive-rank semi-Lagrangian algorithm has been developed that avoid dimensional splitting. 

A third foundational area of relevance to our work is adaptive cross approximation (ACA) for matrices, which employs a greedy sampling strategy to approximate matrices through CUR decompositions \cite{bebendorf2000approximation}.  The CUR decomposition of an \(m \times n\) matrix \(A\) is given by 
    \begin{equation}\label{eq:CUR_first_glance}
    A = \mathcal{CUR},
    \end{equation} 
    with \(\mathcal{C}\) an \(m \times k\) matrix consisting of \(k\) columns of \(A\), \(\mathcal{R}\) a \(k \times n\) matrix containing \(k\) rows of \(A\), and \(\mathcal{U}\) is a \(k \times k\) matrix computed from \(\mathcal{C}\) and \(\mathcal{R}\). 
This decomposition strategy was first explored in \cite{tyrtyshnikov1995pseudo} with an emphasis on efficient selection mechanisms for constructing low-rank approximations. Subsequent advancements included the development of selection strategies, such as maxvol \cite{goreinov1997theory}, discrete empirical
interpolation method (DEIM) \cite{sorensen2016deim,chaturantabut2010nonlinear}, and leverage scores \cite{mahoney2009cur,drineas2008relative}. 
Maxvol will use the skeleton matrix which corresponds to the submatrix with largest modulus determinant, however this is an NP-Hard problem \cite{civril2007finding}. DEIM will utilize information of left and right singular vectors to make robust choices of rows and columns; 
while leverage scores leads to probabilistic strategies to sampling of rows and columns was also studied.
However, these require knowledge of singular vectors, which hinders efficiency of the method for high order tensor decompositions due to curse of dimensionality. 
Along with sampling/selection procedures, there are different options in constructing factorization 
includes projection based components that are shown to be optimal under the Frobenius norm in \cite{stewart1999four}, the recursive update \cite{shi2024distributed} from column and row selection, 
as well as those in \cite{engquist2007fast,cortinovis2024sublinear} which perform projection based on only partial access to the full data. 
Furthermore, CUR decomposition has been extended to the high order tensor setting \cite{grasedyck2015parallel,grasedyck2010hierarchical,shi2024distributed,dolgov2020parallel}. 

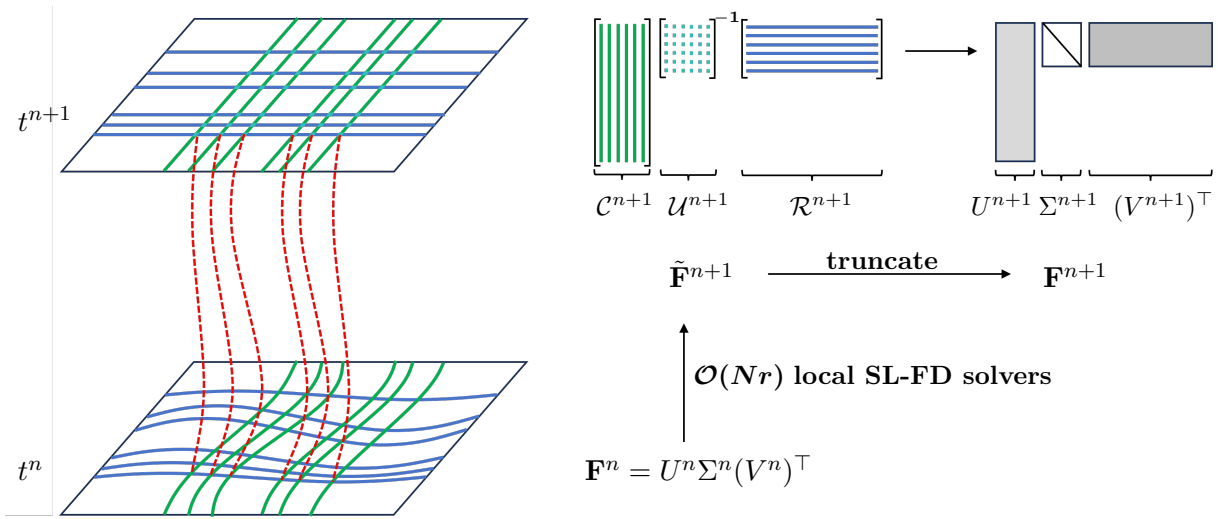
\begin{figure}[h!]
    \begin{center}
        \begin{overpic}[scale=0.59]{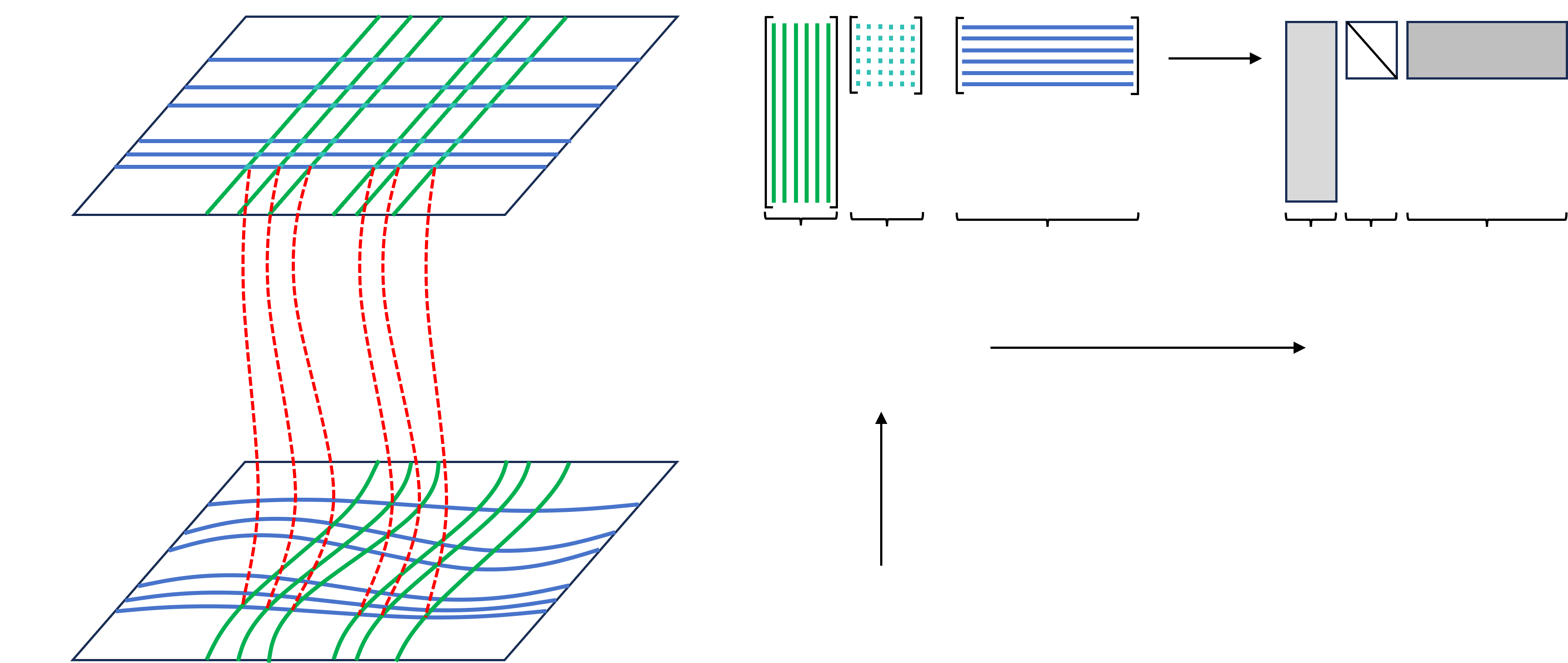}
            \put(1,32){$t^{n+1}$}
            \put(1,3) {$t^n$}
            \put(58.7,41) {\tiny$\mathbf{-1}$}
            \put(49,25) {\small$\mathcal{C}^{n+1}$}
            \put(55,25) {\small$\mathcal{U}^{n+1}$}
            \put(65,25) {\small$\mathcal{R}^{n+1}$}
            \put(80,25) {\small$U^{n+1}$}
            \put(85.6,25) {\small$\Sigma^{n+1}$}
            \put(91.8,25) {\small$(V^{n+1})^\top$}
            \put(48,3) {$\mathbf{F}^n = U^n\Sigma^n(V^n)^\top$}
            \put(55,19) {$\tilde{\mathbf{F}}^{n+1}$}
            \put(86,19) {$\mathbf{F}^{n+1}$}
            \put(68,20.5) {\small\textbf{truncate}}
            \put(57,11) {\small$\boldsymbol{\mathcal{O}(Nr)}$ \textbf{local SL-FD solvers}}
        \end{overpic}
\end{center}
\caption{The SLAR method for linear advection equations.}
\label{fig:SLAR_advection}
\end{figure}

Building on these advancements, we propose a semi-Lagrangian adaptive-rank (SLAR) method without dimensional splitting for linear advection equations and the nonlinear Vlasov-Poisson (VP) system. Our method is built upon a deterministic, non-splitting SL FD solver that achieves up to third-order accuracy.
We assume the solution (on a tensor product of grid points) admits a low rank decomposition $U\Sigma V^\top$, {where $U$ and $V$ are orthonormal singular vectors with decaying singular values as diagonal entries of the diagonal matrix $\Sigma$.} The proposed SLAR method fits into the SAT framework for dynamic rank evolution of the solution, yet our approach distinguishes itself in a 2D setting by employing a matrix sampling strategy of CUR type with effective greedy row and column selections; the sampled row and column elements are updated based on a semi-Lagrangian scheme. 
\Cref{fig:SLAR_advection} schematically illustrate the proposed SLAR procedure. On the left plot, the local SL FD solver is depicted by dashed red curves; this solver operates by tracing information along characteristic curves and applying local polynomial reconstructions at the feet of characteristics. A greedy row/column selection guided by forward-characteristic tracking is used to iteratively select and construct an adaptive cross approximation (ACA) of the solution matrix, which pinpoints which grid points require updates. Note that selected solution matrix entries are updated only on an a as-needed basis, hence the reduced computational complexity (see blue and green lines at the updated time $t^{n+1}$ and the corresponding curves representing the track-back points at the current time step $t^n$). With these updates, the method then constructs the ACA approximation of the solution $C^{n+1}U^{n+1}R^{n+1}$, which is further truncated using singular value decomposition (SVD) to ensure numerical stability. Our method extends to the nonlinear VP system with high-order nonlinear characteristic tracing and spatial interpolation, and we achieve local conservation of charge density by incorporating a correction inspired by the Local Macroscopic Conservative (LoMaC) procedure \cite{guo2024conservative, guo2022local}.

We highlight main novelties of our proposed approach, addressing several existing key challenges in the area. 
\begin{itemize}
    \item {\em High order SL-FD methods without dimensional splitting.}   
Our proposed approach is built upon state-of-art high order SL FD methods in tracking characteristics without dimensional splitting, and with a compact third-order local reconstruction to solution at the characteristic feet. In the nonlinear setting, we employ a Runge-Kutta exponential integrator \cite{cai2018high} for an accurate tracking of characteristics for the nonlinear VP system. Note that preservation of local mass conservation 
in the SL FD framework is challenging. In \cite{xiong2019conservative}, a mass-conservative, non-splitting SL-FD method was developed based on a flux-difference formulation; however, this method introduces an additional time step constraint for numerical stability \cite{xiong2019conservative}. In this work, we 
leverage an implicit scheme for the macroscopic charge density equation, and perform a correction, to enable conservation without additional time step restriction. 

\item {\em Step-and-Truncate via Sampling.} 
Our proposed SLAR algorithm combines the ACA for the 'step' and SVD for the 'truncate' phases. This approach leverages a sampling-enabled strategy to advance the solution over time through characteristic tracing, introducing a unique perspective distinct from Eulerian-type methods, such as those in \cite{kormann2015semi, GuoVlasovFlowMap2022}. Unlike these methods, which require a low-rank decomposition of solutions or related right-hand side (RHS) terms (possibly nonlinear terms) derived from PDE operators, our sampling framework bypasses this requirement. Instead, as long as there is an efficient way to evaluate functions at specific grid points (i.e., matrix/tensor entries) on an as-needed basis, low-rank decomposition is not required. In the SVD 'truncate' step, we propose using a larger truncation threshold, leading to effectively an oversampling algorithm closly related to that of \cite{anderson2015spectral}. This approach offers several advantages, including filtering out high-frequency modes (potentially resulting from numerical interpolation errors) and enhancing the stability of the time-stepping algorithm.
    \item {\em LoMaC correction for the nonlinear VP system.}   
    In the nonlinear VP setting, to ensure local conservation of charge density, we develop an implicit, high-order, and conservative solver for the charge density equation obtained by taking the zeroth moment of the VP system. The implicit discretization is designed to accommodate the larger time steps permitted by the SL scheme. We solve the linear systems arising from the implicit scheme using the Generalized Minimal Residual (GMRES) method with an incomplete LU preconditioner. With this conservative update of the charge density, we apply a LoMaC correction. Due to the unique feature of the sampling-based low rank approach, where low-rank decomposition of RHS terms is unnecessary as long as efficient function evaluation is possible, we propose to perform the LoMaC correction using a Maxwellian defined by local density, momentum, and temperature. This enables a more flexible and robust LoMaC correction compared to the LoMaC scheme proposed for Eulerian methods \cite{guo2022local}.
\end{itemize}

To the best of our knowledge, the proposed SLAR method is the first adaptive-rank SL method, with high order (up to third order) accuracy in both spatial and time, allowing for large time stepping size, without dimensional splitting, and with preservation of local conservation laws. In this paper, we focus exclusively on the 2-D matrix case to exploit the  algorithm design.
The rest of the paper is organized as follows. \Cref{sec:EL_FV_WENO_method} presents the proposed SLAR method with subsections discussing linear and nonlinear problems, elaborating technical details of the proposed procedure; \Cref{sec:numerical_tests} showcases a variety of numerical tests demonstrating the effectiveness of the proposed SL methods; finally we conclude in \Cref{sec:conclusion}.

\section{SLAR Methods}\label{sec:EL_FV_WENO_method}

In this section, we first propose the SLAR method for linear advection equations in \Cref{sec:SLAR_LT}; followed by introducing 
the LoMaC SLAR method for the nonlinear VP system in \Cref{sec:SLAR_VP}.

\subsection{SLAR method for linear advection equations}\label{sec:SLAR_LT}
Consider 
\begin{equation}\label{eq:linear_advection_2d}
    f_t + a(x,y,t)f_x + b(x,y,t)f_y = 0,~~x\in[x_L,x_R],~~y\in[y_B,y_T]
\end{equation}
with $(a(x,y,t),b(x,y,t))$ being a known velocity field. 
We assume uniform meshes for the $x$- and $y$-dimensions 
$$x_L = x_{\frac{1}{2}} < x_{\frac{3}{2}} < \ldots < x_{N_x+\frac{1}{2}} = x_R,$$ 
$$y_B = y_{\frac{1}{2}} < y_{\frac{3}{2}} < \ldots < y_{N_y+\frac{1}{2}} = y_T,$$ 
with $x_i = \frac{1}{2}(x_{i-\frac{1}{2}}+x_{i+\frac{1}{2}})$, $\Delta x = x_{i+\frac{1}{2}} - x_{i-\frac{1}{2}}$, $I^x_i = [x_{i-\frac{1}{2}},x_{i+\frac{1}{2}}]$, $y_j = \frac{1}{2}(y_{j-\frac{1}{2}}+y_{j+\frac{1}{2}})$, $\Delta y = y_{j+\frac{1}{2}} - y_{j-\frac{1}{2}}$, $I^y_j = [y_{j-\frac{1}{2}},y_{j+\frac{1}{2}}]$, and $I_{i,j} = I^x_i\times I^y_j$ for all $i,j$. We also let $\mathbf{x} := [x_1, x_2,\ldots, x_{N_x}]^\top$ and $\mathbf{y} := [y_1, y_2,\ldots, y_{N_y}]^\top$.

Throughout our scheme design, we assume that the initial condition of \eqref{eq:linear_advection_2d}, as well as the solution, can be approximated by the following low-rank format:
\begin{equation}\label{eq:low_rank_ini_assum1}
    f(x,y,t) = \sum_{k=1}^r \sigma_k(t) u_k(x, t)v_k(y, t),
\end{equation}
where \(\{u_k(x, t)\}_{k=1}^r\) and \(\{v_k(y, t)\}_{k=1}^r\) are time-dependent orthonormal basis functions in their respective dimensions, and \(\{\sigma_k(t)\}_{k=1}^r\) are time-dependent coefficients. This low-rank form of \eqref{eq:low_rank_ini_assum1} corresponds to an SVD of the solution matrix living on the 2D tensor product grid, i.e.,
\begin{equation}\label{eq:low_rank_ini_assum2}
    \mathbf{F}^n = U^n\Sigma^n (V^n)^{\top} = \sum_{k=1}^r \sigma^n_k \mathbf{u}^n_k(\mathbf{v}^n_k)^{\top},
\end{equation}
where $\mathbf{F}^n\in \mathbb{R}^{N_x\times N_y}$ represents the solution, with superscript $n$ corresponding to the solution snapshot at time $t^n$ with explicit time dependence on basis vectors as well as coefficients. We will skip such superscript for notational brevity, when there is no ambiguity. Here, $U = [\mathbf{u}_1, \ldots, \mathbf{u}_{r}]\in \mathbb{R}^{N_x\times r}$, $\Sigma = \text{diag}\{\sigma_1, \ldots, \sigma_{r}\}\in \mathbb{R}^{r\times r}$, $V = [\mathbf{v}_1, \ldots , \mathbf{v}_{r}]\in \mathbb{R}^{N_y\times r}$. When $r\ll \min\{N_x, N_y\}$, such representation offers significant compression in data storage and computational efficiency. Eq.~\eqref{eq:low_rank_ini_assum2} is the basic form of solution we start from, and return to, in each time step evolution. In the following, we will elaborate three key components of the SLAR algorithm for linear problems: the SL FD update in \Cref{sec:SLFD_method}, the ACA of updated solution matrix in \Cref{sec:Adaptive_cross}, and the SVD truncation for numerical stability and overall algorithm summary in \Cref{sec:AR_Cross_SVD}.

\subsubsection{Local non-splitting SL-FD method}\label{sec:SLFD_method}
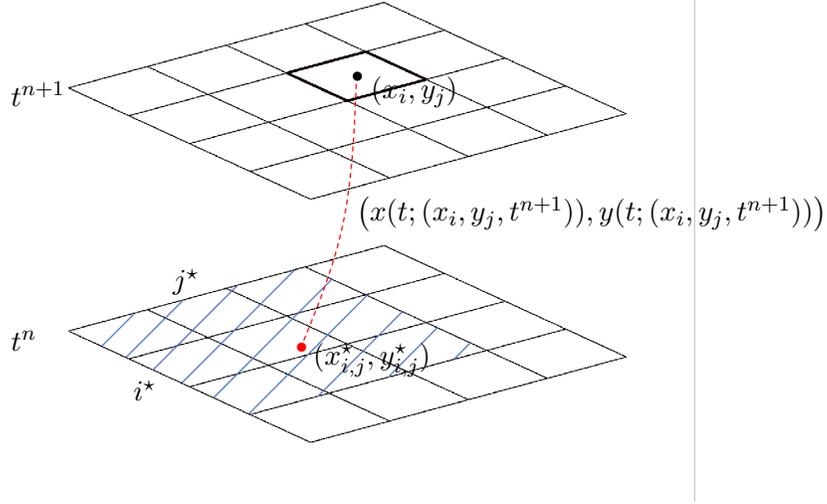
\begin{figure}[htb]
	\begin{center}
		\begin{overpic}[scale=0.5]{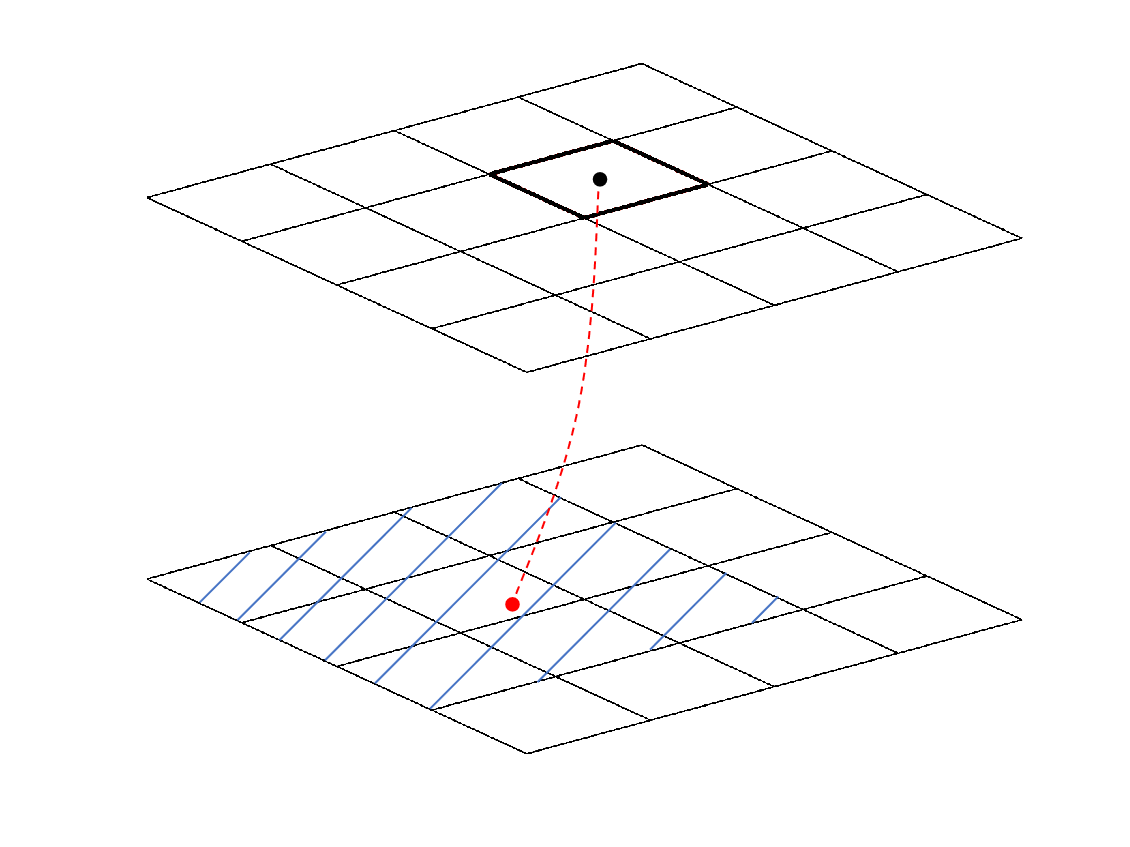}
			\put(55,56){$(x_i,y_j)$}
			\put(5,55){$t^{n+1}$}
			\put(5,21) {$t^n$}
			\put(47,19){$(x^\star_{i, j}, y^\star_{i, j})$}
                \put(22,14) {$i^\star$}
                \put(27.5,29.5) {$j^\star$}
                \put(53,39) {$\left(x(t; (x_i, y_j,t^{n+1})), y(t; (x_i, y_j,t^{n+1}))\right)$}
		\end{overpic}
	\end{center}
	\caption{Schematic illustration of tracing characteristics.}
	\label{fig:characteristic_tracing}
\end{figure}

Below, we present details of the non-splitting SL FD method, which is used to evaluate local matrix entries during the sampling step. It is well known that the solution of \eqref{eq:linear_advection_2d} can be obtained by tracing characteristics backward:
\begin{equation}\label{eq:SL_FD_nondiscretiztion}
    f(x_i,y_j,t^{n+1}) = f(x^\star_{i,j}, y^\star_{i,j},t^n),
\end{equation}
where $(x^\star_{i, j}, y^\star_{i, j}) = \left(x(t^n; (x_i, y_j, t^{n+1})), y(t^n; (x_i, y_j, t^{n+1}))\right)$ is the characteristic foot at $t^n$. Here, $\left(x(t; (x_i, y_j, t^{n+1})), y(t; (x_i, y_j, t^{n+1}))\right)$ represents the characteristic curve passing through $(x_i, y_j, t^{n+1})$ (see \Cref{fig:characteristic_tracing}), satisfying the following system of equations:
\begin{equation}\label{eq:characteristic_ODE}
    \begin{cases}
        \frac{d(x(t))}{dt} = a(x(t), y(t), t),\\
        \frac{d(y(t))}{dt} = b(x(t), y(t), t),\\
        x(t^{n+1}) = x_i,\\
        y(t^{n+1}) = y_j,
    \end{cases}
\end{equation}
which can be solved using a high-order Runge-Kutta (RK) method. 
With \eqref{eq:SL_FD_nondiscretiztion}, a simple SL-FD method can be implemented as follows:
\begin{equation}\label{eq:full_rank_SL}
    f_{i,j}^{n+1} = \mathcal{R}\left(\{f^n_{p,q}\}_{(p,q) \in \mathcal{I}(i^\star,j^\star)}\right) (x_{i,j}^\star, y_{i,j}^\star),
\end{equation}
where $\mathcal{R}$ is a local third-order reconstruction operator that depends on 9 local nodal values, $(i^\star,j^\star)$ is the index pair such that $(x_{i,j}^\star, y_{i,j}^\star) \in I_{i^\star, j^\star}$, and $\mathcal{I}(i^\star,j^\star):=\{(i^\star+l_1, j^\star+l_2)\}_{l_1, l_2=-1}^{1}$ defines a $3 \times 3$ stencil centered around $(i^\star, j^\star)$ (see the shaded area in \Cref{fig:characteristic_tracing}).

To construct a robust and efficient third-order reconstruction operator, we define
\begin{equation}
    V := \{p(x, y) \in P^2(I_{i^\star,j^\star}) | p(x_k, y_l) = f^n_{k,l}, \quad \text{for}~ (k,l) \in \mathcal{S}_{\text{interpolation}}\},
\end{equation}
where $\mathcal{S}_{\text{interpolation}} := \{(i^\star, j^\star), (i^\star-1, j^\star), (i^\star+1, j^\star), (i^\star, j^\star-1), (i^\star, j^\star+1)\}$. The reconstruction operator is then given by
\begin{equation}
    \mathcal{R}\left(\{f^n_{p,q}\}_{(p,q) \in \mathcal{I}(i^\star,j^\star)}\right)(x, y) = \min_{p \in V}\left[\sum_{(k,l) \in \mathcal{S}_{\text{least\_square}}}\left(p(x_k, y_l) - f^n_{k,l}\right)^2\right]^{\frac{1}{2}},
\end{equation}
where $\mathcal{S}_{\text{least\_square}} := \mathcal{I}(i^\star,j^\star) \setminus \mathcal{S}_{\text{interpolation}}$. The least-square procedure yields the explicit polynomial $p(x, y) = \sum_{l=1}^6 a_l P_l(x, y)$, where the coefficients $a_l$ are given by:
\begin{equation}
    \begin{split}
        a_1 &= f^n_{i,j}, \quad a_2 = \frac{1}{2}(f^n_{i+1,j} - f^n_{i-1,j}), \quad a_3 = \frac{1}{2}(f^n_{i,j+1} - f^n_{i,j-1}), \\
        a_4 &= -f^n_{i,j} + \frac{1}{2}(f^n_{i-1,j} + f^n_{i+1,j}), \\
        a_5 &= \frac{1}{4}(f^n_{i+1,j+1} + f^n_{i-1,j-1} - f^n_{i-1,j+1} - f^n_{i+1,j-1}), \\
        a_6 &= -f^n_{i,j} + \frac{1}{2}(f^n_{i,j-1} + f^n_{i,j+1}),
    \end{split}
\end{equation}
with the local polynomial basis:
\begin{equation}
\begin{split}
    P_1(x,y) &= 1, \quad P_2(x,y) = \xi_i(x), \quad P_3(x,y) = \eta_j(y), \\
    P_4(x,y) &= \xi_i(x)^2, \quad P_5(x,y) = \xi_i(x) \eta_j(y), \quad P_6(x,y) = \eta_j(y)^2,
\end{split}
\end{equation}
where $\xi_i(x) = \frac{x - x_i}{\Delta x}$ and $\eta_j(y) = \frac{y - y_j}{\Delta y}$. Finally, we remark that the SL FD algorithm has issues of mass conservation in a general nonlinear setting \cite{xiong2019conservative}, which we will address in later in the paper.  
 
\subsubsection{ACA of matrices}\label{sec:Adaptive_cross}
\begin{figure}[ht]
    \centering
    \includegraphics[width=0.65\linewidth]{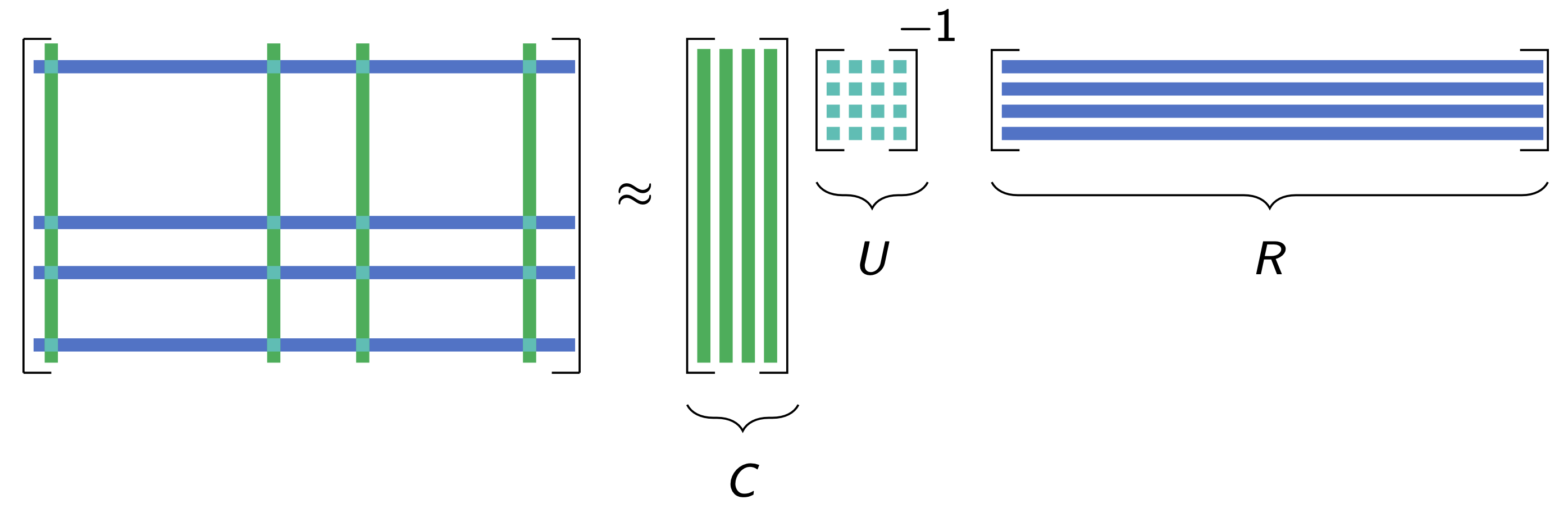}
    \caption{Visual representation of the CUR decomposition used in SLAR. }
    \label{fig:CUR_Picture}
\end{figure}

Next, we provide a complete description of the ACA, which is one of the key components of the full SLAR method. The ACA algorithm is based on a CUR decomposition of an $m\times n$ matrix $A$ in the form of \eqref{eq:CUR_first_glance}. In a standard CUR decomposition, $\mathcal{C}$ and $\mathcal{R}$ are selected columns and rows of $A$, with row and column indices denoted by $\mathcal{I}$ and $\mathcal{J}$; $\mathcal{U}$ is the inverse of $A(\mathcal{I}, \mathcal{J})$ (i.e. the intersection of rows and columns of $A$). Figure~\ref{fig:CUR_Picture} illustrates the CUR decomposition, with selected rows and columns highlighted in blue and green, and $A(\mathcal{I}, \mathcal{J})$ highlighted in turquoise. The selected rows and columns, together with local update from the SL scheme, enable an efficient algorithm that only need to access a small percent of the full data, thereby reducing the computing time and storage requirements of the proposed SLAR. For implementation of CUR, working with $\mathcal{U}$, directly as $A(\mathcal{I}, \mathcal{J})^{-1}$, may not be effective, as it may introduce significant numerical errors due to large condition numbers of the intersection matrix.


We propose to use an adaptive recursive algorithm to realize the ACA of matrices. The algorithm is summarized in \Cref{alg:ar_cross} with a prescribed threshold of $\varepsilon_C$. 
Assume we have a rank $k-1$ approximation of a matrix $A$, denoted as $A_{k-1}$, as well as row and column index sets $\mathcal{I}, \mathcal{J}$, in the k-th iteration step of ACA. There are two phases in updating $A_k$ as well as the index sets: 
\begin{itemize}
    \item {\bf Phase I: row and column selection.} We use a greedy pivot selection. This algorithm will search for the entry with the largest residual referred to as a pivot. In order to start the algorithm we will make a simple sampling step in which we select $p$ random points in line 3 of~\Cref{alg:ar_cross} that do not coincide with any row or column that has already been selected. It is a greedy algorithm, as we can't state that the selection is the largest residual, but rather an entry that has largest residual in its respective row and column. This can be seen in lines 4 and 5 of~\Cref{alg:ar_cross}. 
    Once a new row and column index is identified, we then expand the index sets, see line 6 of ~\Cref{alg:ar_cross}. 

\item {\bf Phase II: construction of $\mathcal{CUR}$.} In SLAR, we use a recursive rank-1 update, as in Proposition~\ref{prop:rank_one_update} \cite{shi2024distributed}, as an efficient and robust update of $\mathcal{CUR}$, as we make selections for $\mathcal{C}$ and $\mathcal{R}$.
In the implementation, we use \Cref{cor:E_JDE_I} to explicitly construct the ACA of updated solution matrix. A crucial consequence of~\Cref{cor:E_JDE_I} is that we can express the CUR decomposition in a form that mimics the SVD in the sense that we have a tall skinny matrix $\mathcal{E}_\mathcal{J}$ times a small diagonal matrix $\mathcal{D}$ times a short fat matrix $\mathcal{E}_\mathcal{I}$. 
\end{itemize}


\begin{algorithm}[h]
\caption{ACA algorithm}
\label{alg:ar_cross}
\begin{algorithmic}[1]
\Statex \textbf{Input:} Access to entries of $A\in N_x\times N_y$ (direct or function); initial row index range $\mathbb{I}_x$; initial column index range $\mathbb{I}_y$; tolerance $\varepsilon_{\text{C}}$; maximum rank $r_{\text{max}}$.
\Statex \textbf{Output:} $\tilde{A} = \mathcal{E}_\mathcal{J}\mathcal{D}\mathcal{E}_\mathcal{I}$; $\mathcal{I}$ and $\mathcal{J}$.
\State Set $A_0 = [0]_{N_x\times N_y}$,  $\mathcal{I} = \emptyset$ and $\mathcal{J} = \emptyset$
\State \textbf{for} $k = 1, 2, \ldots, r_{\text{max}}$
\State\quad\quad Pick $p$ samples $\mathcal{L} = \{(i_l,j_l)\}_{l=1}^p$  randomly with $i_l\in\mathbb{I}_x\setminus \mathcal{I}$ and $j_l\in\mathbb{I}_y\setminus\mathcal{J}$ and choose the one with largest error:
$$(i^*, j^*) \leftarrow \argmax_{(i, j) \in \mathcal{L}} |A_{i,j} - (A_{k-1})_{i,j}|$$
\State\quad\quad Search a column for maximum error:
$$i_k^* \leftarrow \argmax_{i \in \mathbb{I}_x\setminus\mathcal{I}} |A_{i,j^*} - (A_{k-1})_{i,j^*}|$$
\State\quad\quad Search a row for maximum error:
$$j_k^* \leftarrow \argmax_{j \in \mathbb{I}_y\setminus\mathcal{J}} |A_{i_k^*,j} - (A_{k-1})_{i_k^*,j}|$$
\State\quad\quad Expand the index sets: $\mathcal{I} \leftarrow \mathcal{I}\cup\{i_k^*\};\quad \mathcal{J} \leftarrow \mathcal{J}\cup\{j_k^*\}$ 
\State\quad\quad Update rank-$k$ approximation with~\Cref{eq:Recursive_CUR}.

\State\quad\quad \textbf{If} $\|\frac{1}{(E_{k-1})_{i_k^\star,j_k^\star}}(E_{k-1})_{:,j_k^\star}(E_{k-1})_{i_k^\star,:}\|_F < \varepsilon_{\text{C}} $ or $|(E_{k-1})_{i_k^\star,j_k^\star}| < 10^{-14}$
\State\quad\quad\quad\quad \textbf{break}
\State\quad\quad \textbf{End}
\State \textbf{End}
\State $\tilde{A} \leftarrow A_k$
\end{algorithmic}
\end{algorithm}

\begin{prop}[Recursive update of cross \cite{shi2024distributed}]\label{prop:rank_one_update}
 Assume we have a rank - $(k-1)$ cross approximation $A_{k-1} = C_{k-1}U_{k-1}R_{k-1}$ for row and column indices $\mathcal{I}$ and $\mathcal{J}$, then the cross approximation $A_k$ for rows $\mathcal{I}\cup \{i\}$ and columns $\mathcal{J}\cup \{j\}$ is given by
\begin{equation}\label{eq:Recursive_CUR}
A_{k} = A_{k-1} + \frac{1}{(A-A_{k-1})_{i,j}}(A-A_{k-1})_{:,j}(A-A_{k-1})_{i,:}.
\end{equation}
\end{prop}
\begin{cor}\label{cor:E_JDE_I}
For $\mathcal{I} = \{i_1,i_2,\dots,i_k\}, \mathcal{J} = \{j_1,j_2,\dots,j_k\}$, ~\cref{eq:Recursive_CUR} can be equivalently represented as
$$
A_k = \mathcal{E}_{\mathcal{J}}D\mathcal{E}_\mathcal{I}
$$
\begin{align*}
\mathcal{E}_\mathcal{J} & = [e_{0,:,j_1},e_{1,:,j_2},\dots,e_{k-1,:,j_{k}}] \in \mathbb{R}^{N_x\times r}\\
D &= \text{diag}\left(e_{0,i_1,j_1}^{-1},e_{1,i_2,j_2}^{-1},\dots,e_{k-1,i_{k},j_{k}}^{-1}\right)\in \mathbb{R}^{r\times r}\\
\mathcal{E}_\mathcal{I} &= [e_{0,i_1,:}^T,e_{1,i_2,:}^T,\dots,e_{k-1,i_{k},:}^T]^T\in \mathbb{R}^{r\times N_y}
\end{align*}
with the notation $e_{l,m,n} = A(m,n) - A_{l}(m,n)$ and ``$:$" signifies all entries.
\end{cor}

In the context of solving time-dependent PDEs, one could choose to narrow down index ranges $\mathbb{I}_x$ and $\mathbb{I}_y$ for the pivot search.  
For example, for linear advection equation with solution that has a compact support (see examples in the numerical section), we could dynamically estimate index ranges $\mathbb{I}_x$ and $\mathbb{I}_y$ using a forward Lagrangian characteristic tracing. 
In particular, we define a rectangular region $\Omega^n := [x_{i^n_m},x_{i^n_M}]\times[y_{j^n_m},y_{j^n_M}]$, and the associated index ranges \(\mathbb{I}^n_x = \{i^n_m, i^n_m+1, \ldots, i^n_M\}\) and \(\mathbb{I}^n_y = \{j^n_m, j^n_m+1, \ldots, j^n_M\}\), such that
\begin{equation}
 \begin{split}
     &i^{n}_m \leftarrow \max\{1, \min_{i\in \mathcal{I}} \{i\}-2\};\quad~ i^{n}_M \leftarrow \min\{N_x, \max_{i\in \mathcal{I}} \{i\}+2\};\\
     &j^{n}_m \leftarrow \max\{1, \min_{j\in \mathcal{J}} \{j\}-2\};\quad j^{n}_M \leftarrow \min\{N_y, \max_{j\in \mathcal{J}} \{j\}+2\},
 \end{split}
 \end{equation}
where \(\mathcal{I}\) and \(\mathcal{J}\) are the selected row and column index sets at \(t^n\), respectively, and \(\pm 2\) is used to slightly expand the ranges in case they are too narrow. The ranges, \(\mathbb{I}^n_x\) and \(\mathbb{I}^n_y\), approximate the index ranges for which the numerical solution is non-zero. We then perform the forward Lagrangian characteristic tracing procedure, summarized in \Cref{alg:predict_new_index_ranges}, to predict narrowed index ranges $\tilde{\mathbb{I}}^{n+1}_x$ and $\tilde{\mathbb{I}}^{n+1}_y$ for the pivot search at \(t^{n+1}\).


\begin{algorithm}[h]
\caption{Algorithm of predicting new index ranges}
\label{alg:predict_new_index_ranges}
\begin{algorithmic}[1]
\Statex \textbf{Input:} Current row index range $\mathbb{I}^n_x$; current column index range $\mathbb{I}^n_y$; current time $t^n$; future time $t^{n+1}$; velocity field $(a(x,y,t),b(x,y,t))$; x-dimension mesh $\mathbf{x}$; y-dimension mesh $\mathbf{y}$.
\Statex \textbf{Output:} Predicted new index ranges $\tilde{\mathbb{I}}^{n+1}_x$ and $\tilde{\mathbb{I}}^{n+1}_y$ at $t^{n+1}$.
\State Randomly sample $s$ points on each edge of the $\partial \Omega^n$ and form a sample set:
$$\mathcal{B}^n \leftarrow \{(x^n_k,y^n_k)\}_{k=1}^{4s}$$
\State Extend the sample set $\mathcal{B}^n$ by adding the four vertices of $\partial \Omega^n$:
$$\mathcal{B}^n \leftarrow \mathcal{B}^n\cup\{(x_{i^n_m},y_{j^n_m}),(x_{i_m^n},y_{j^n_M}),(x_{i^n_M},y_{j^n_m}),(x_{i^n_M},y_{j^n_M})\}$$
\State Trace the characteristic curves that cross the points in $\mathcal{B}^n$ forward and find the feet set $\tilde{\mathcal{B}}^{n+1}$ at $t^{n+1}$ by solving:
\begin{equation*}
\begin{cases}
    dX(t)/dt = a(X(t),Y(t),t),\\
    dY(t)/dt = b(X(t),Y(t),t),\\
    \left(X(t^n),Y(t^n)\right) \in \mathcal{B}^n
\end{cases}
\end{equation*}
using the third-order RK method
\State Predict new index ranges $\tilde{\mathbb{I}}^{n+1}_x=(\tilde{i}^{n+1}_m,\tilde{i}^{n+1}_m+1,\ldots,\tilde{i}^{n+1}_M)$ and $\tilde{\mathbb{I}}^{n+1}_y = (\tilde{j}^{n+1}_m,\tilde{j}^{n+1}_m+1,\ldots,\tilde{j}^{n+1}_M)$:
\begin{equation*}
\small
\begin{split}
    &\tilde{i}^{n+1}_m \leftarrow \max\{i\in\{1,\ldots,N_x\}|x_{i} \leq \min_{(x,y)\in\tilde{\mathcal{B}}^{n+1}}\{x\}\},\quad \tilde{i}^{n+1}_M \leftarrow \min\{i\in\{1,\ldots,N_x\}|x_{i} \geq \max_{(x,y)\in\tilde{\mathcal{B}}^{n+1}}\{x\}\},\\
    &\tilde{j}^{n+1}_m \leftarrow \max\{j\in\{1,\ldots,N_y\}|y_{j} \leq \min_{(x,y)\in\tilde{\mathcal{B}}^{n+1}}\{y\}\},\quad \tilde{j}^{n+1}_M \leftarrow \min\{j\in\{1,\ldots,N_y\}|y_{j} \geq \max_{(x,y)\in\tilde{\mathcal{B}}^{n+1}}\{y\}\}
\end{split}
\end{equation*}
\end{algorithmic}
\end{algorithm}

\subsubsection{SLAR algorithm with SVD truncation}\label{sec:AR_Cross_SVD}

The SL FD update of matrix entries and the ACA of matrices just described provide basic ingredients in the proposed SLAR algorithm. To enhance the stability, we propose to perform an additional SVD truncation step, with the truncation threshold $\varepsilon_S$ larger than $\varepsilon_C$ for the CUR decomposition. 

We summarize the SLAR method in \Cref{alg:SLAR}. In our description, we denote the local SL-FD method as 
\begin{equation}
    f^{n+1}_{i,j} = \text{SL}_{i,j}(\mathbf{F}^n,t^n,t^{n+1},a(x,y,t),b(x,y,t)),
\end{equation}
where $\mathbf{F}^{n} = U^n\Sigma^n(V^n)^\top\in\mathbb{R}^{N_x\times N_y}$ is the SVD solution at $t^n$. 
The proposed method uses the ACA algorithm to construct the cross approximation 
\[
\tilde{\mathbf{F}}^{n+1} = \mathcal{E}_{\mathcal{J}}^{n+1}D^{n+1}\mathcal{E}^{n+1}_\mathcal{I}
\]
which uses 
\[
\text{SL}_{\cdot,\cdot}\left(\mathbf{F}^n,t^n,t^{n+1},(a(x,y,t),b(x,y,t)\right)
\]
to access updated solution at arbitrary grid points in an as-needed basis, as specified in line 1 of \Cref{alg:SLAR}.
 Following this, an efficient SVD truncation is performed to stabilize the SLAR method. The SVD truncation involves applying QR factorization on $\mathcal{E}^{n+1}_{\mathcal{J}}$ and $\mathcal{E}^{n+1}_{\mathcal{I}}$, followed by applying a standard SVD to the small $r\times r$ matrix $D^{n+1}$, as detailed from line 2 to line 5 in \Cref{alg:SLAR}. The resulting $\mathbf{F}^{n+1} = U^{n+1}\Sigma^{n+1}(V^{n+1})^{\top}$ is the truncated SVD decomposition of updated solution
 with tolerance $\epsilon_S$. The SVD truncation effectively truncate modes with smaller singular values, some of which may be caused from numerical discretization error or artificial oscillations. The SVD truncation reduces the Frobenius norm, enhancing stability as an approximation for PDE solutions. 
 In practice, we require $\varepsilon_C < \varepsilon_S$, where $\varepsilon_C$ and $\varepsilon_S$ are the tolerances for the ACA and the SVD truncation, respectively. This naturally leads to the result $r_{S} < r_{C}$. 

In summary, the SLAR update in \Cref{alg:SLAR} is denoted by
\begin{equation}\label{eq:SLAR}
    \mathbf{F}^{n+1} = \text{SLAR}(\mathbf{F}^n,t^n,t^{n+1},a(x,y,t),b(x,y,t),\varepsilon_C,r_{\max},\varepsilon_S).
\end{equation}
with $\mathbf{F}$ represented in a low rank form as $U \Sigma V^\top$.

 




\begin{algorithm}[h]
\caption{SLAR method}
\label{alg:SLAR}
\begin{algorithmic}[1]
\Statex \textbf{Input:} Current SVD solution $\mathbf{F}^n = U^n\Sigma^n (V^n)^\top$; current time $t^n$; future time $t^{n+1}$; velocity field $(a(x,y,t),b(x,y,t))$; tolerance $\varepsilon_{\text{C}}$ of the ACA; maximum rank $r_{\text{max}}$ of the ACA; tolerance $\varepsilon_S$ of the SVD truncation.
\Statex \textbf{Output:} SVD solution $\mathbf{F}^{n+1} = U^{n+1}\Sigma^{n+1}(V^{n+1})^{\top}$.
\State Use $\text{SL}_{\cdot,\cdot}\left(\mathbf{F}^n,t^n,t^{n+1},(a(x,y,t),b(x,y,t)\right)$, 
$\mathbb{I}_x$, $\mathbb{I}_y$, $\varepsilon_{\text{C}}$, and $r_{\text{max}}$ as input, and \textbf{call \Cref{alg:ar_cross}} to update a cross approximation $\tilde{\mathbf{F}}^{n+1} = \mathcal{E}_{\mathcal{J}}^{n+1}D^{n+1}\mathcal{E}^{n+1}_\mathcal{I}$ 
\State Apply QR factorization to $\mathcal{E}^{n+1}_{\mathcal{J}}$ and $(\mathcal{E}^{n+1}_{\mathcal{I}})^{\top}$:
\begin{equation*}
    \mathcal{E}^{n+1}_{\mathcal{J}} = Q_1R_1,~~(\mathcal{E}^{n+1}_{\mathcal{I}})^{\top} = Q_2R_2
\end{equation*}
\State Apply SVD to $R_1D^{n+1}R_2^{\top}$:
\begin{equation*}
    R_1D^{n+1}R_2^{\top} = U_S\Sigma_S V_S^{\top}
\end{equation*}
\State Determine the rank for tolerance $\varepsilon_S$:
$$r_S \leftarrow \min\{k=1,2,\ldots,r_C|(\Sigma_S)_{k+1,k+1} > \varepsilon_{S}\}$$
\State Construct the SVD solution $\mathbf{F}^{n+1} = U^{n+1}\Sigma^{n+1}(V^{n+1})^{\top}$:
\begin{equation*}
    U^{n+1} \leftarrow Q_1(U_S)_{:,1:r_S},~~\Sigma^{n+1}\leftarrow(\Sigma_S)_{1:r_S,1:r_S},~~V^{n+1}\leftarrow Q_2(V_S)_{:,1:r_S}
\end{equation*}
\end{algorithmic}
\end{algorithm}

\subsection{LoMaC SLAR method for the nonlinear VP system}\label{sec:SLAR_VP}
In this subsection, we generalize SLAR to a nonlinear 1D1V VP system, 
\begin{equation}\label{eq:VP_1}
f_t + vf_x + E(x,t)f_v = 0,\quad x\in\Omega_x,\quad v\in\mathbb{R},
\end{equation}
\begin{equation}\label{eq:VP_2}
E(x,t) = -\phi_x,\quad -\phi_{xx}(x,t) = \rho(x,t) - \rho_0,
\end{equation}
where \((x,v)\) is the coordinate of the phase space, \(f(x,v,t)\) is the probability distribution function of finding a particle at position \(x\) with velocity \(v\) at time \(t\), \(E\) is the electric field, \(\phi\) is the electrostatic potential,  \(\rho=\int_{\mathbb{R}}f(x,v,0)\,dv\) is the charge density, and \(\rho_0 = \frac{1}{\Omega_x}\int_{\Omega_x}\int_{\mathbb{R}}f(x,v,0)\,dv\,dx\). 
Similar discretization with $(x_i, v_j)$ as a given grid point is used. 

The SLAR method for the nonlinear VP system has two technical issues to address: one on nonlinear characteristic tracing, see discussions in \Cref{sec:SLAR_VP_CF3}; another one on ensuring local conservation of charge density, as will be discussed in \Cref{sec:lomac}.




\subsubsection{RK exponential integrators for nonlinear characteristics tracing}\label{sec:SLAR_VP_CF3}
To track characteristics for the nonlinear VP system, we apply the Runge Kutta exponential integrators with up to third order temporal accuracy. 
Such RK expontential integrators freeze the velocity field at RK stages, for which a linear solver can be directly used, please see \cite{cai2021high} for detailed discussion. For example, the following Butcher table, associated with a third order RK exponential integrator, 
\begin{equation}
\setlength\arraycolsep{10pt} 
\renewcommand{\arraystretch}{1.5} 
\begin{array}{c|ccc}
0 & & & \\
\frac{1}{3} & \frac{1}{3}  & & \\
\frac{2}{3} & 0 & \frac{2}{3} & \\
\hline
 & \frac{1}{3}  & 0 & 0 \\
 & -\frac{1}{12} & 0 &\frac{3}{4}
\end{array}
\end{equation}
decompose the nonlinear characteristics tracing into solving a sequence of linear advection equations (with frozen velocity field taken as a linear combination of fields from previous RK stages) in the following fashion.
\begin{equation}\label{eq:CF3}
\begin{split}
    \mathbf{F}^{(1)} &= \text{SLAR}(\mathbf{F}^n,t^n,t^{n+1},\frac{1}{3}v,\frac{1}{3}E^n(x),\varepsilon_C,r_{\max},\varepsilon_S),\\
    \mathbf{F}^{(2)} &= \text{SLAR}(\mathbf{F}^n,t^n,t^{n+1},\frac{2}{3}v,\frac{2}{3}E^{(1)}(x),\varepsilon_C,r_{\max},\varepsilon_S),\\
    \mathbf{F}^{n+1,\star} &= \text{SLAR}(\mathbf{F}^{(1)},t^n,t^{n+1},\frac{2}{3}v,-\frac{1}{12}E^{n}(x) + \frac{3}{4}E^{(2)}(x),\varepsilon_C,r_{\max},\varepsilon_S),
\end{split}
\end{equation}
where \(E^n(x)\), \(E^{(1)}(x)\), and \(E^{(2)}(x)\) are obtained from \(\mathbf{F}^n\), \(\mathbf{F}^{(1)}\), and \(\mathbf{F}^{(2)}\) using a Poisson solver.  
The \(\star\) symbol in \(\mathbf{F}^{n+1,\star}\) indicates that this is not yet the final updated solution; a LoMaC type adjustment will be in place to ensure local conservation of charge density. 

Taking the first three moments of the SVD distribution \(F^{n+1,\star} = U\Sigma V^\top\), we estimate the discrete macroscopic charge, current, and kinetic energy densities, \(\boldsymbol{\rho}^{n+1,\star}\), \(\mathbf{J}^{n+1,\star}\), and \(\boldsymbol{\kappa}^{n+1,\star} \in \mathbb{R}^{N_x}\), at the updated time level, by
\begin{equation}
\begin{cases}
        \boldsymbol{\rho}^{n+1,\star}= \Delta v U\Sigma V^\top \mathbf{1}_v,\\
        \mathbf{J}^{n+1,\star} =  \Delta v U\Sigma V^\top \mathbf{v},\\
        \boldsymbol{\kappa}^{n+1,\star} =  \Delta v U\Sigma V^\top (\frac{1}{2}\mathbf{v}^2),\\
    \end{cases}
    \label{eq: macro}
\end{equation}
where \(\mathbf{1}_v \in \mathbb{R}^{N_v}\) is the vector of all ones, \(\mathbf{v}\in\mathbb{R}^{N_v}\) is the vector of grid points in \(v\), and \(\mathbf{v}^2\) is the element-wise square of \(\mathbf{v}\).
These macroscopic quantities will be used, for linearization and prediction, in the next subsection to design a LoMaC correction.

\subsubsection{LoMaC SLAR method}\label{sec:lomac}

In this subsubsection, we propose to implicitly evolve the charge density equation of the VP system (as zeroth moment of the VP system), 
\begin{equation}\label{eq:charge_density_eq}
    \rho_t + (\rho u)_x = 0,
\end{equation}
where \(\rho u = {\int_{\mathbb{R}} f v\, dv}\). 
The implicit nature of the scheme, allows for a large time stepping size, aligned with that for the SLAR algorithm. The scheme is obtained by first approximating \((\rho u)_x\) with an upwind discretization via flux splitting, followed by a third-order, stiffly accurate, diagonally implicit Runge-Kutta (DIRK) method for the time integration. Then a LoMaC type correction \cite{guo2022local} will be in place to correct the low rank solution in order to ensure local conservation of charge density. 

The flux splitting for the flux $g \doteq \rho u$ is designed as follows: \(g = g^+ + g^-\), where \(g^\pm = \frac{u \pm \alpha}{2}\rho\) with \(\alpha = \max\{|u|\}\).  We define upwind differential matrices as follows:
\begin{equation}\label{eq:diff_matrix_1} 
   \mathbf{D}^{\mathbf{u}}\boldsymbol{\rho} := \mathbf{D^+}\,\left(\text{diag}\left\{\frac{\mathbf{u} + \alpha}{2}\right\} \boldsymbol{\rho}\right) + \mathbf{D^-}\, \left(\text{diag}\left\{\frac{\mathbf{u} - \alpha}{2}\right\}\boldsymbol{\rho}\right) =  \mathbf{D^+}\mathbf{g}^{+} + \mathbf{D^-}\mathbf{g}^{-},
\end{equation}
where \(\mathbf{D}^+\) is assembled using the following left-biased interpolation for the positive splitting of velocity
\begin{equation}
    {(\mathbf{D}^+ \mathbf{g}^{+})_{i} :=}  \frac{1}{\Delta x}\left[\frac{1}{6}g^+_{i-2} - g^+_{i-1} + \frac{1}{2}g^+_i + \frac{1}{3}g^+_{i+1}\right],
\end{equation}
and \(\mathbf{D}^-\) is assembled for that of the negative splitting 
\begin{equation}
   {(\mathbf{D}^- \mathbf{g}^{-})_i :=} \frac{1}{\Delta x}\left[-\frac{1}{3}g^-_{i-1} - \frac{1}{2}g^-_i + g^-_{i+1} - \frac{1}{6}g^-_{i+2}\right].
\end{equation}

To illustrate the idea of implicit solver, we consider a first order backward Euler as a prototype scheme. 
\[
\boldsymbol{\rho}^{n+1} = \boldsymbol{\rho}^n - \Delta t\mathbf{D}^{\mathbf{u}^{n+1,\star}}\boldsymbol{\rho}^{n+1}
\]
where {$\mathbf{u}^{n+1,\star} = \frac{\mathbf{J}^{n+1,\star}}{\boldsymbol{\rho}^{n+1,\star}}$} is the predicted macroscopic velocity obtained from eq.~\eqref{eq: macro} explicitly.  
To solve such an implicit scheme, we apply the sparse GMRES method with an incomplete LU preconditioner with a drop tolerance of \(10^{-6}\) and initial guesses \(\boldsymbol{\rho}^{n+1,\star}\) from eq.~\eqref{eq: macro} \cite{elman1982iterative,saad1986gmres}.

To extend to high order time discretization, we use the thrid-order, stiffly accurate DIRK method with the following Butcher table \cite{alexander1977diagonally}:
\begin{equation}
\label{eq:S-stable_DIRK3}
\setlength\arraycolsep{10pt} 
\renewcommand{\arraystretch}{1.5} 
\begin{array}{c|ccc}
 \beta & \beta & 0 & 0\\
 \tau_2 & \tau_2 - \beta & \beta & 0 \\
 1 & b_1 & b_2 & \beta\\
 \hline
 & b_1 & b_2 & \beta
\end{array}
\end{equation}
where \(\beta\) is the root of \(x^3 - 3x^2 + \frac{3}{2}x - \frac{1}{6} = 0\) lying in \((\frac{1}{6}, \frac{1}{2})\), \(\tau_2 = \frac{1 + \beta}{2}\), \(b_1 = -\frac{6\beta^2 - 16\beta + 1}{4}\), and \(b_2 = \frac{6\beta^2 - 20\beta + 5}{4}\). To predict the macroscopic $\boldsymbol{\rho}$ and $\mathbf{u}$ for initial guesses and for a linearized velocity field, we perform an interpolation at intermediate RK stages from corresponding quantities at $t^{n-1}$, $t^n$ and predicted ones at $t^{n+1}$ from eq.~\eqref{eq: macro}.

To enable local conservation of charge density, we adjust the non-conservative kinetic solution \(\mathbf{F}^{n+1,\star}\), using the conservative charge density \(\boldsymbol{\rho}^{n+1}\) computed from the implicit update of \eqref{eq:charge_density_eq} as proposed above. We define two local Maxwellians, 
$\mathbf{M}^{n+1}, \mathbf{M}^{n+1,\star} \in \mathbf{R}^{N_x\times N_v}$
with $(i, j)$-entry being
\begin{equation}
    M^{n+1}_{i,j} = \frac{\rho^{n+1}_i}{2\pi T^{n+1,\star}_i}\exp\left(\frac{-|v_j - u^{n+1,\star}_i|^2}{2T^{n+1,\star}_i}\right),~~\text{and}~~M^{n+1,\star}_{i,j} = \frac{\rho^{n+1,\star}_i}{2\pi T^{n+1,\star}_i}\exp\left(\frac{-|v_j - u^{n+1,\star}_i|^2}{2T^{n+1,\star}_i}\right) \nonumber
\end{equation}
where \(\mathbf{T}^{n+1,\star}=(T^{n+1,\star}_i) \in \mathbb{R}^{N_x}\) is the predicted temperature at \(t^{n+1}\), with \(T^{n+1,\star}_i = 2\kappa^{n+1,\star}_i/\rho^{n+1,\star}_i - (u^{n+1,\star}_i)^2\).
 Finally, we conduct the LoMaC correction to the SVD distribution by
\begin{equation}\label{eq:conservative_SLAR}
\begin{split}
       \mathbf{F}^{n+1} &= \mathbf{F}^{n+1,\star} + \mathbf{M}^{n+1} - \mathbf{M}^{n+1,\star},\\
      &= \mathbf{F}^{n+1,\star} + \left[\left(\boldsymbol{\rho}^{n+1} - \boldsymbol{\rho}^{n+1,\star}\right)./\boldsymbol{\rho}^{n+1}\right].*\mathbf{M}^{n+1},
\end{split}
\end{equation}
where ``$./$" and ``$.*$" denote element-wise division and multiplication, respectively. 
The final correction term in \eqref{eq:conservative_SLAR} can be interpreted as to adjust the local charge density with a localized Maxwellian distribution function.
The numerical solution at the end of a time step update of the LoMaC SLAR method is a summation of low rank prediction and an explicit correction term in Maxwellian form (no necessarily in the low rank format). 

\section{Numerical tests}\label{sec:numerical_tests}
In this section, we present the benchmark results for linear advection equations in \Cref{sec:numer_linear}. The numerical tests for the VP system are provided in \Cref{sec:numer_VP}. Unless otherwise specified, we use the same tolerance settings of \(\varepsilon_C = 10^{-4}\), \(\varepsilon_S = 10^{-3}\), and no maximum rank limitation is applied for the ACA algorithm. The time step size is determined by 
\begin{equation}
    \Delta t = \frac{\text{CFL}}{\left(\frac{\max\{|a|\}}{\Delta x} + \frac{\max\{|b|\}}{\Delta y}\right)},
\end{equation}
where \(\max\{|a|\}\) and \(\max\{|b|\}\) represent the exact maximum absolute values of the velocity field for linear advection simulations, and the maximum absolute values of the discrete velocity on the spatial grid for the Vlasov Poisson system. All boundary conditions are either periodic or zero-boundary. The exploration of various practical boundary conditions will be left as a topic for future work in our ongoing research.

\subsection{Linear advection equations}\label{sec:numer_linear}
In this subsection, we present three benchmark tests: the linear advection equation with constant coefficients, the rigid body rotation, and the swirling deformation flow. Through these tests, we aim to investigate the spatial and temporal order of accuracy, the adaptive-rank behavior, and the compression ratio of the degrees of freedom (DOFs), defined as the ratio of the total entries in the SVD solution to the total entries in the full matrix solution.

\begin{exa}(2-D advection equation with constant coefficients).  Consider the equation
\begin{equation}\label{2d_linear_advection}
  u_t + u_x + u_y = 0,\quad x,~y\in[-\pi,\pi],
\end{equation}
with the initial condition \(u(x,y,0) = \sin(x+y)\). The exact solution for this problem is \(u(x,y,t) = \sin(x+y-2t)\). We set the final time to \(T=2\) and use a CFL number of 1. We provide the \(L^1\) and \(L^{\infty}\) errors for varying meshes, along with the corresponding orders of accuracy, in \Cref{tab_2d_linear_advection}. 
The designed spatial order is third, so the results in the table reflect the correct order of accuracy. We also present the average ranks of the SVD and the cross approximation for different mesh settings. As shown, the SVD ranks are 2, which is optimal for \(\sin(x+y-2T)\). The cross ranks range from 3 to 5, which may due to numerical errors; the SVD truncation effectively remove redundant modes from numerical errors.

\begin{table}[!htbp]
\centering
\caption{(2-D advection equation with constant coefficients).\label{tab3} $L^1$ and $L^{\infty}$ errors, corresponding orders of accuracy, average SVD ranks, and average CUR ranks at $T = 2$ for $CFL = 1$.}\label{tab_2d_linear_advection}

  \centering
\begin{tabular}{|c|cc|cc|c|c|}
\hline
mesh&$L^1$ error&order&$L^{\infty}$ error&order&SVD rank& cross rank\\
\hline
$8\times8$ & 2.21e-01 & --- & 3.62e-01 & --- & 2.00 & 3.00 \\
\hline
$16 \times 16$& 1.99e-02 & 3.47 & 3.31e-02 & 3.45 & 2.00 & 3.00 \\
\hline
$32 \times 32$& 1.25e-03 & 4.00 & 2.29e-03 & 3.86 & 2.00 & 4.94 \\
\hline
$64 \times 64$& 7.15e-05 & 4.12 & 1.30e-04 & 4.14 & 2.00 & 4.97 \\
\hline
$128 \times 128$& 4.64e-06 & 3.95 & 1.80e-05 & 2.85 & 2.00 & 3.00 \\
\hline
$256 \times 256$& 4.07e-07 & 3.51 & 1.74e-06 & 3.37 & 2.00 & 3.00 \\
\hline
\end{tabular}
\end{table}
\end{exa}

\begin{exa}(Rigid body rotation). Consider the equation
\begin{equation}\label{2d_RBD}
  u_t - (yu)_x + (xu)_y = 0,\quad x,~y\in[-\pi,\pi],
\end{equation}
with the following initial condition,
\begin{equation}\label{sdf_smooth_initial_condtion}
  u(x,y,0) = \begin{cases}
  r^b_0\cos\left(\frac{r^b(\mathbf{x})\pi}{2r^b_0}\right)^6, &\text{if}~ r^b(\mathbf{x})<r^b_0,\\
  0, &\text{otherwise},
  \end{cases}
\end{equation}
where \(r^b_0=0.3\pi\), \(r^b(\mathbf{x})=\sqrt{ (x-x_0^b)^2+(y-y_0^b)^2 }\), and the center of the cosine bell is \((x_0^b,y_0^b) = (0.3\pi,0)\). With a period of \(2\pi\), the cosine bell retains its shape, rotates around \((0,0)\), and returns to its initial position. 

For this test, we set the final time to \(T = 2\pi\), use two fixed meshes (\(128\times 128\) and \(256\times 256\)), and vary the CFL number, which controls the time step size \(\Delta t\). This setup allows us to evaluate the temporal accuracy of the proposed SLAR method. With sufficiently large \(\Delta t\), the error is expected to be dominated by temporal discretization.
As shown in \Cref{fig:CFL_RBR}, the slopes of approximately 3 confirm the expected temporal accuracy. When \(\text{CFL} < 11\), the accumulation of spatial error becomes more significant, dominating the total error. Consequently, we observe that the total error decreases as the time step size increases until the CFL reaches around 11. Under the same settings, we present a semi-log plot of CFL numbers versus average SVD and cross ranks on the left side of \Cref{fig:aver_rank_RBR}, and a semi-log plot of CFL numbers versus average SVD compression ratios of DOFs on the right side. 
We observe that 
the average cross ranks increase with the CFL numbers; while the average SVD rank stays low.  
The extra rank in the ACA algorithm seems to be influenced by 
numerical approximation errors, while the SVD truncation effectively removes redundant modes and noisy modes possibly from numerical errors. On the right side of \Cref{fig:aver_rank_RBR}, the results show that a more refined mesh leads to a better compression ratio of DOFs.

Another interesting setup for the rigid body rotation can be configured by assigning an initial condition \(u(x,y,0) = \exp(-25x^2)\exp(-2y^2)\), see the left plot in \Cref{fig:RBRii_aver_rank}. We use a \(256\times 256\) mesh and a CFL number of 10 to simulate this problem. The rank history of the SVD and cross approximations is presented on the right side of \Cref{fig:RBRii_aver_rank}. As shown, the SVD ranks adaptively increase and decrease according to the orientation of the solution. When \(t=0\), \(\pi/2\), \(\pi\), \(3\pi/2\), and \(2\pi\), the compressed structure is vertical or horizontal, causing the rank to decrease. Conversely, when \(t=\pi/4\), \(3\pi/4\), \(5\pi/4\), and \(7\pi/4\), the compressed structure is diagonal, leading to an increase in rank.

\begin{figure}[!htbp]
\centering
  \subfloat{
  \includegraphics[width=0.4\textwidth]{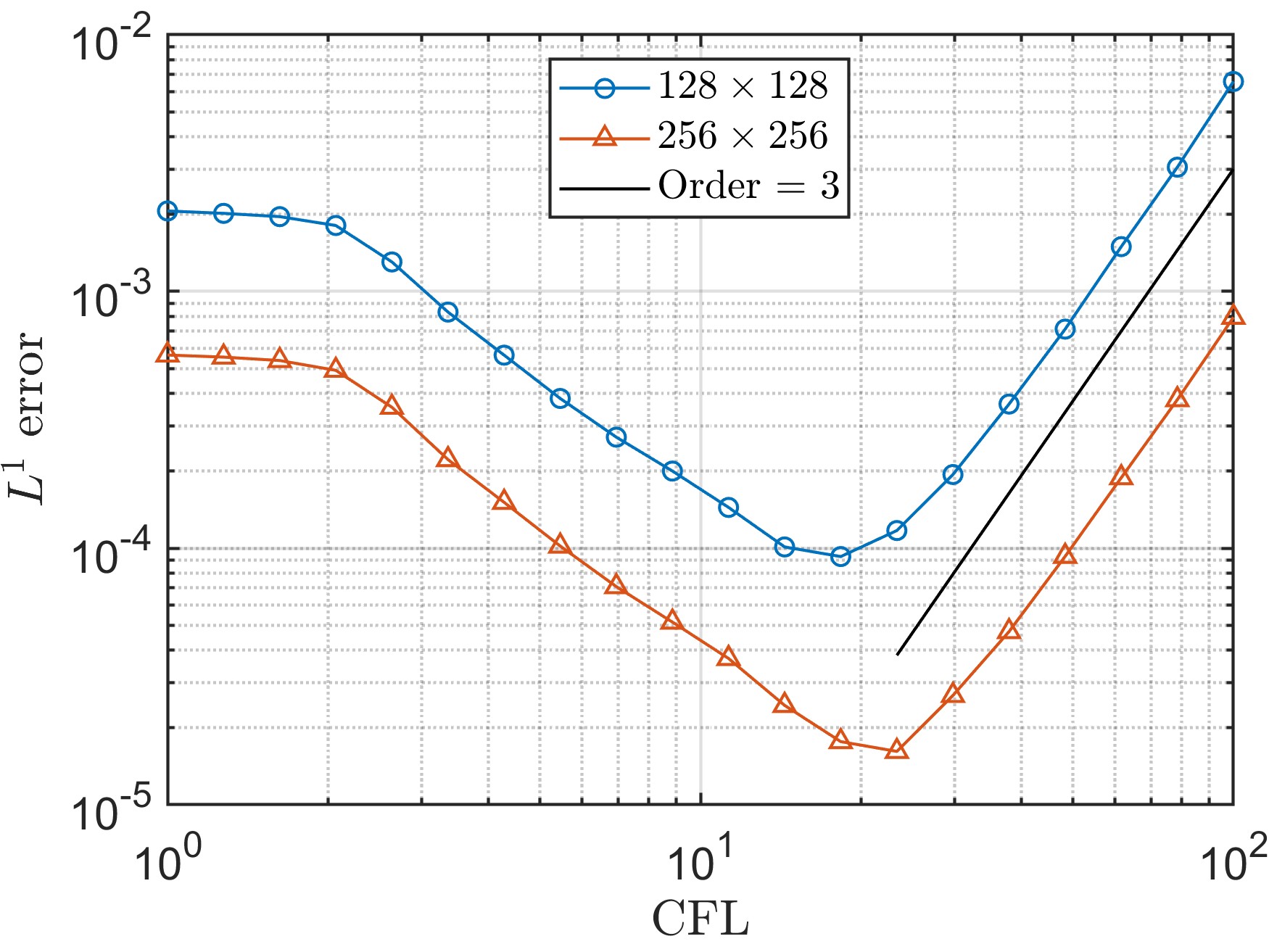}
  }
  \subfloat{
  \includegraphics[width=0.4\textwidth]{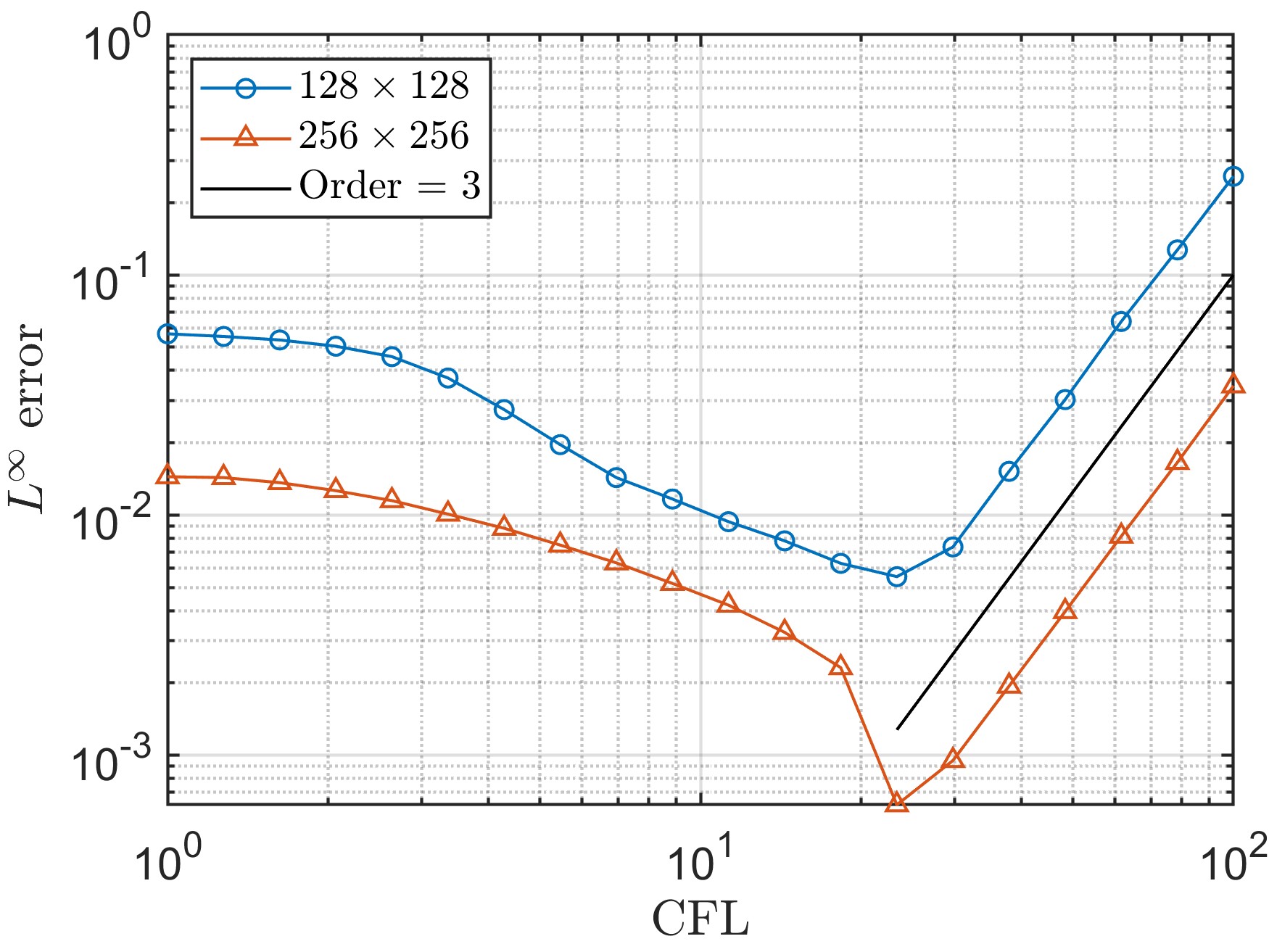}
  }

 \caption{(Rigid body rotation). Log-log plots of CFL numbers versus $L^1$ and $L^{\infty}$ errors with two sets of fixed meshes, $128\times128$ and $256\times256$ at $t = 2\pi$. }
  \label{fig:CFL_RBR}
\end{figure}

\begin{figure}[!htbp]
\centering
  \subfloat{
  \includegraphics[width=0.4\textwidth]{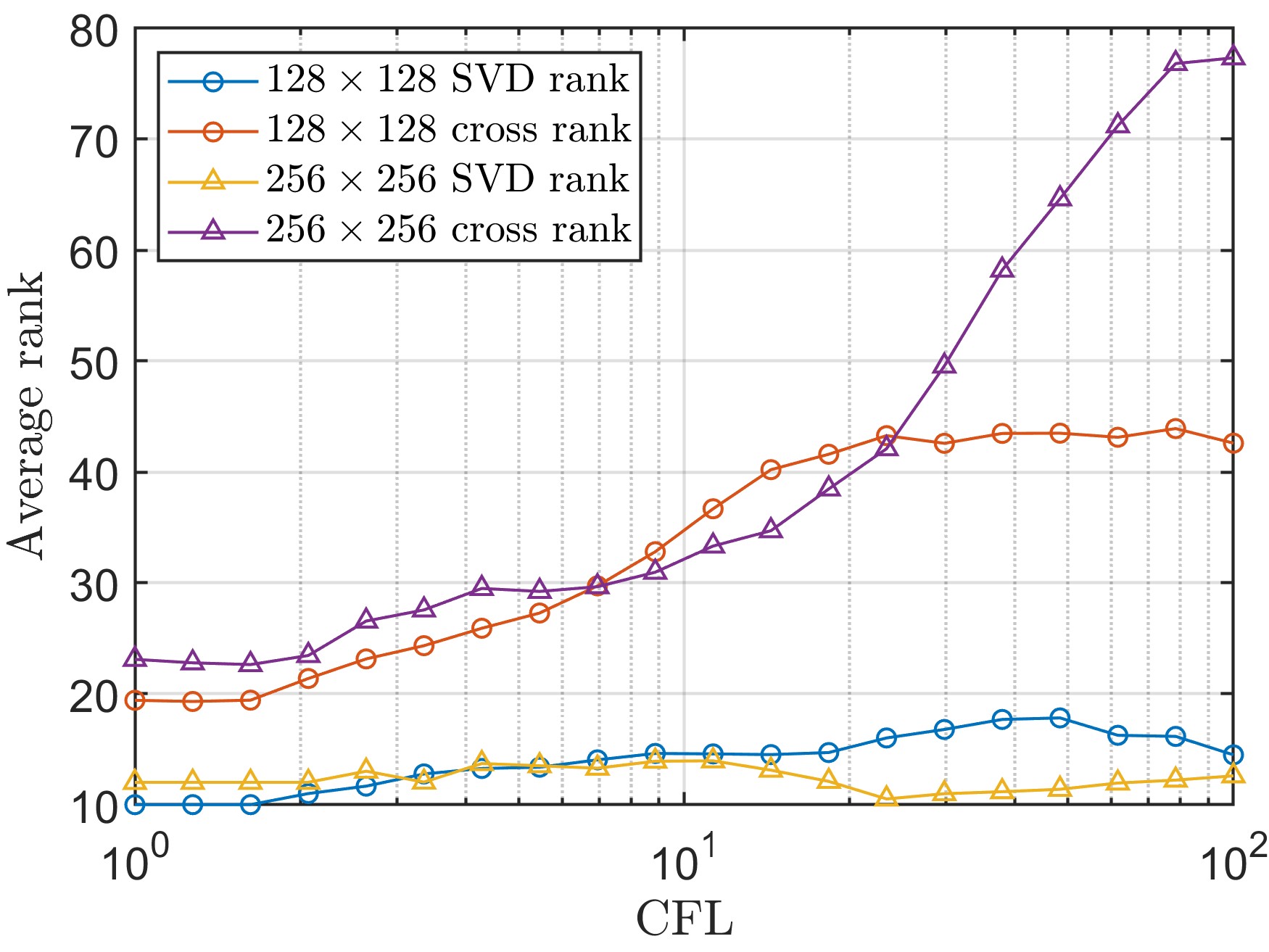}
  }
  \subfloat{
  \includegraphics[width=0.4\textwidth]{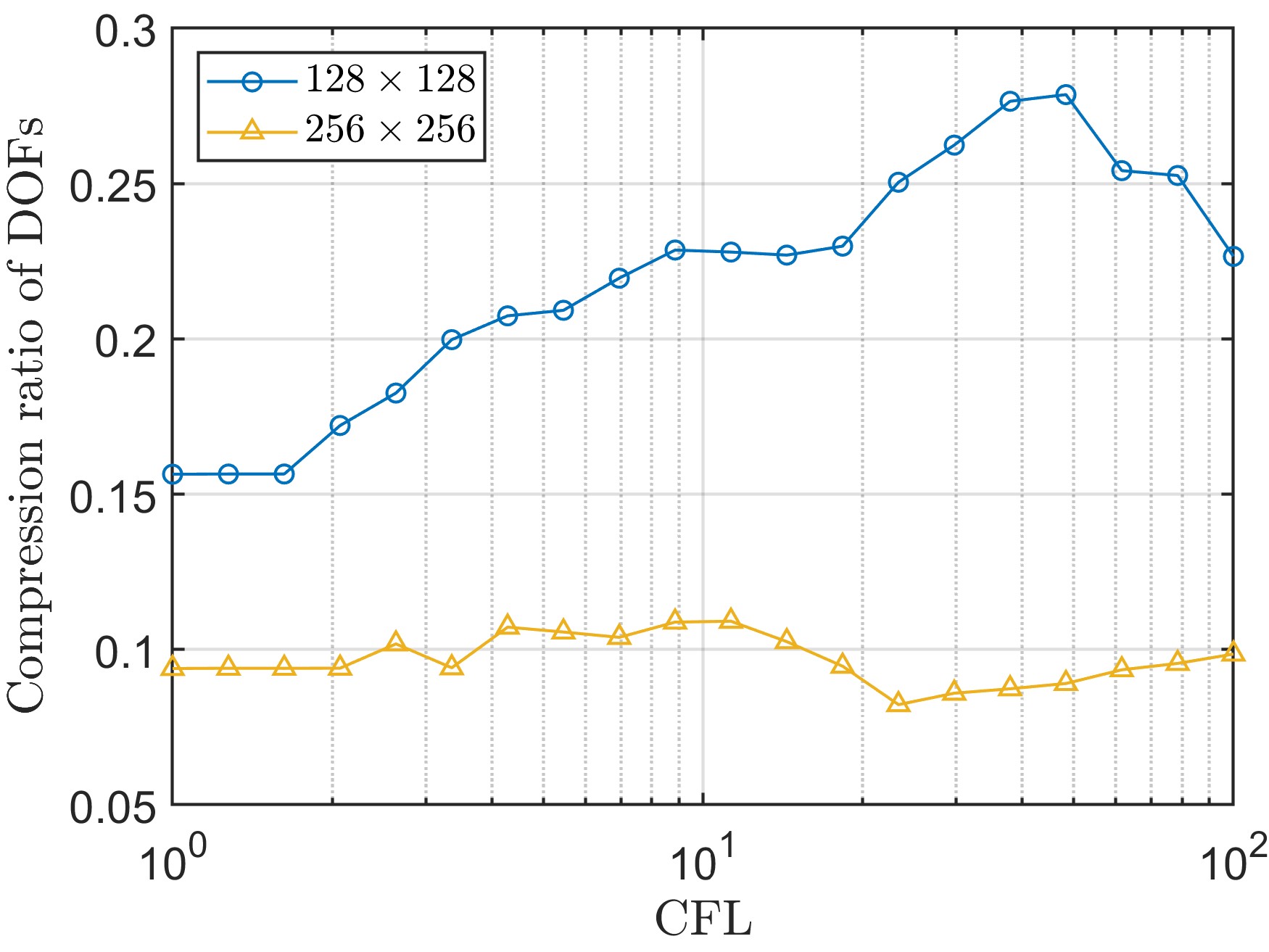}
  }
 \caption{(Rigid body rotation). Left: semi-log plot of CFL numbers versus average ranks of the simulations in \Cref{fig:CFL_RBR}. Right: semi-log plot of CFL numbers versus compression ratios of DOFs of the simulations in \Cref{fig:CFL_RBR}.}
  \label{fig:aver_rank_RBR}
\end{figure}

\begin{figure}[!htbp]
\centering
  \subfloat{
  \includegraphics[width=0.4\textwidth]{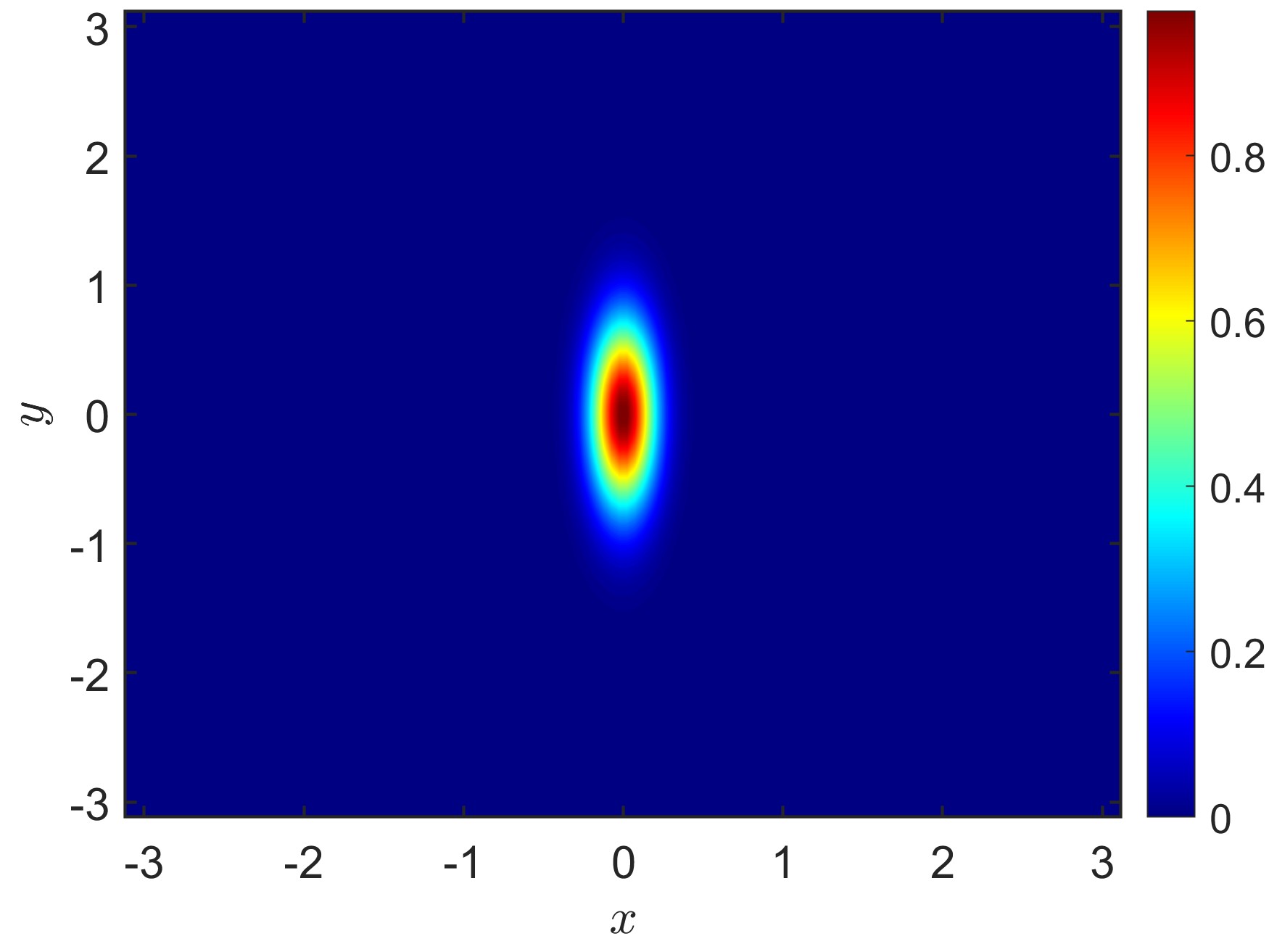}
  }
  \subfloat{
  \includegraphics[width=0.4\textwidth]{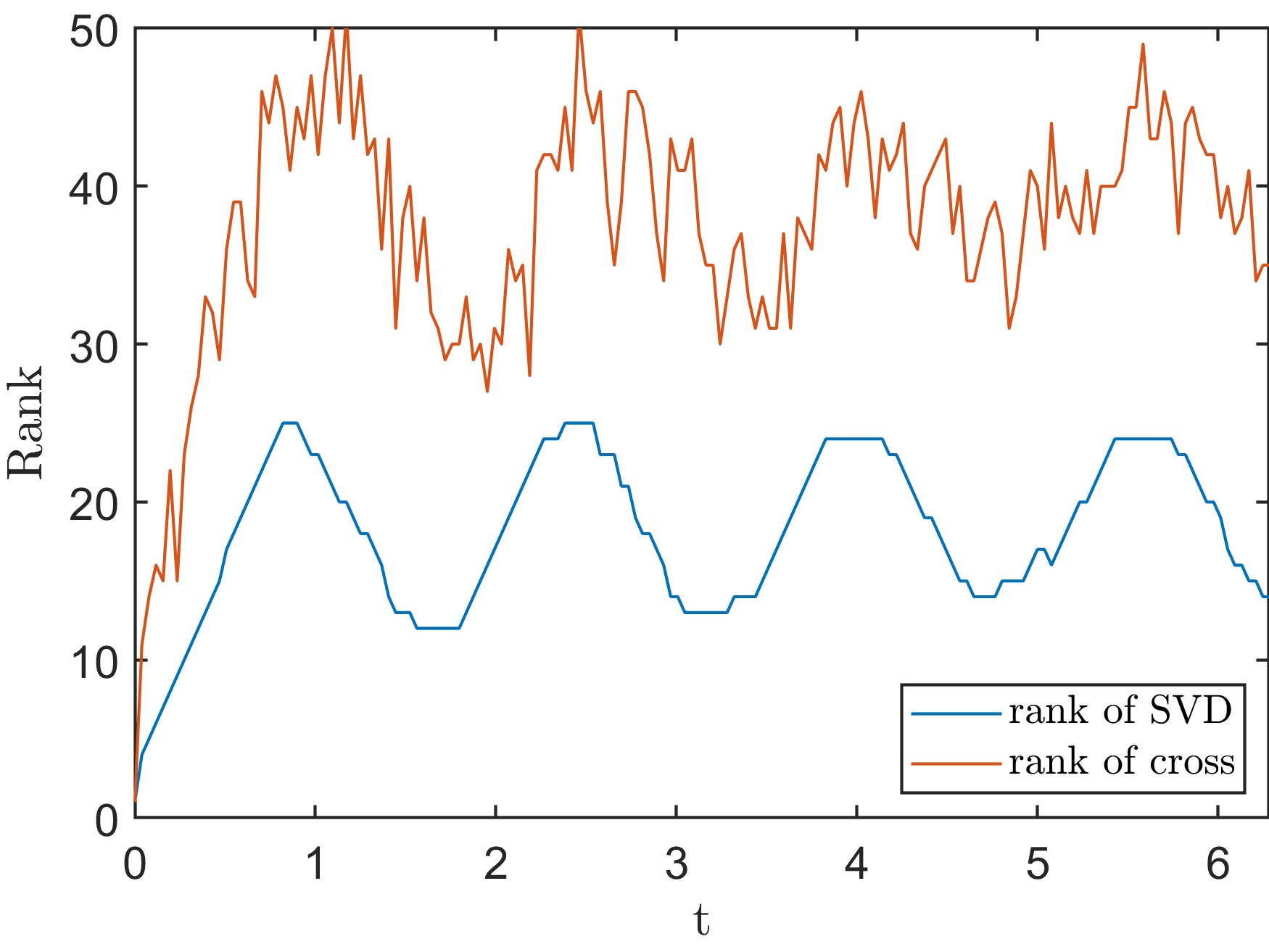}
  }

 \caption{(Rigid body rotation). Left: contour plot of the new initial condition. Right:  rank history of the simulation with a mesh of $256\times256$ and a CFL of 10. }
  \label{fig:RBRii_aver_rank}
\end{figure}
\end{exa}

\begin{exa}(Swirling deformation flow). Consider the equation
\begin{equation}\label{2_D_SDF}
		\begin{split}
			u_t-(2\pi\cos^2\left(\frac x2\right)\sin(y)g(t)u)_x + (2\pi\sin(x)\cos^2\left(\frac y2\right)g(t)u)_y=0,\quad x,~y\in[-\pi,\pi],
		\end{split}
	\end{equation}
where \(g(t) = \cos(\pi t/T)\) with \(T = 1.5\). Note that \eqref{2_D_SDF} is divergence-free, and it remains a linear advection equation. The swirling deformation flow deforms the solution until the half period, \(t = 0.75\), and then reverses the deformation back to the initial condition at the full period, \(t = 1.5\). We use the same initial condition as in \eqref{sdf_smooth_initial_condtion} for this problem, with a \(256\times 256\) mesh and a CFL number of 10.

At the top of \Cref{fig:SDF}, we show the selected columns and rows (left) and the contour plot of the SVD solution (right) at the half period. The contour plot accurately captures the teardrop shape of the solution, while the selected columns and rows illustrate the cross approximation—demonstrating that only a subset of columns and rows is needed to form an effective approximation. On the bottom left of \Cref{fig:SDF}, we visualize the bound prediction procedure from \Cref{alg:predict_new_index_ranges} at the half period of the simulation. The red estimated \(t^{n+1}\) bound is slightly larger and shifted relative to the previous bound following the forward characteristic tracing procedure of \Cref{alg:predict_new_index_ranges}. Lastly, we present the rank history of the simulation for a full period (bottom right). The SVD rank remains relatively stable, while the cross rank fluctuates.

\begin{figure}[h]
\centering
  \subfloat{
  \includegraphics[width=0.4\textwidth]{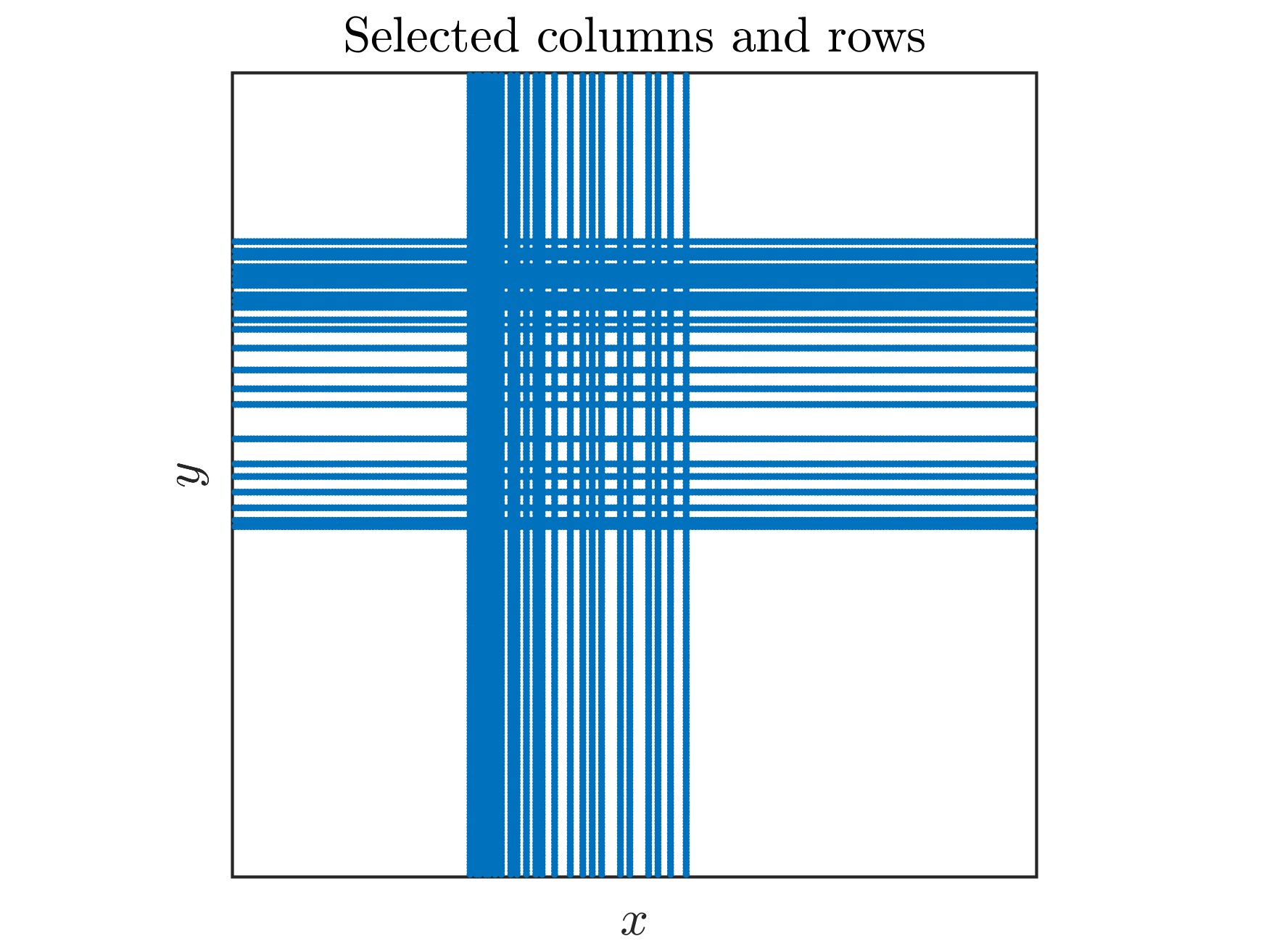}
  }
  \subfloat{
  \includegraphics[width=0.4\textwidth]{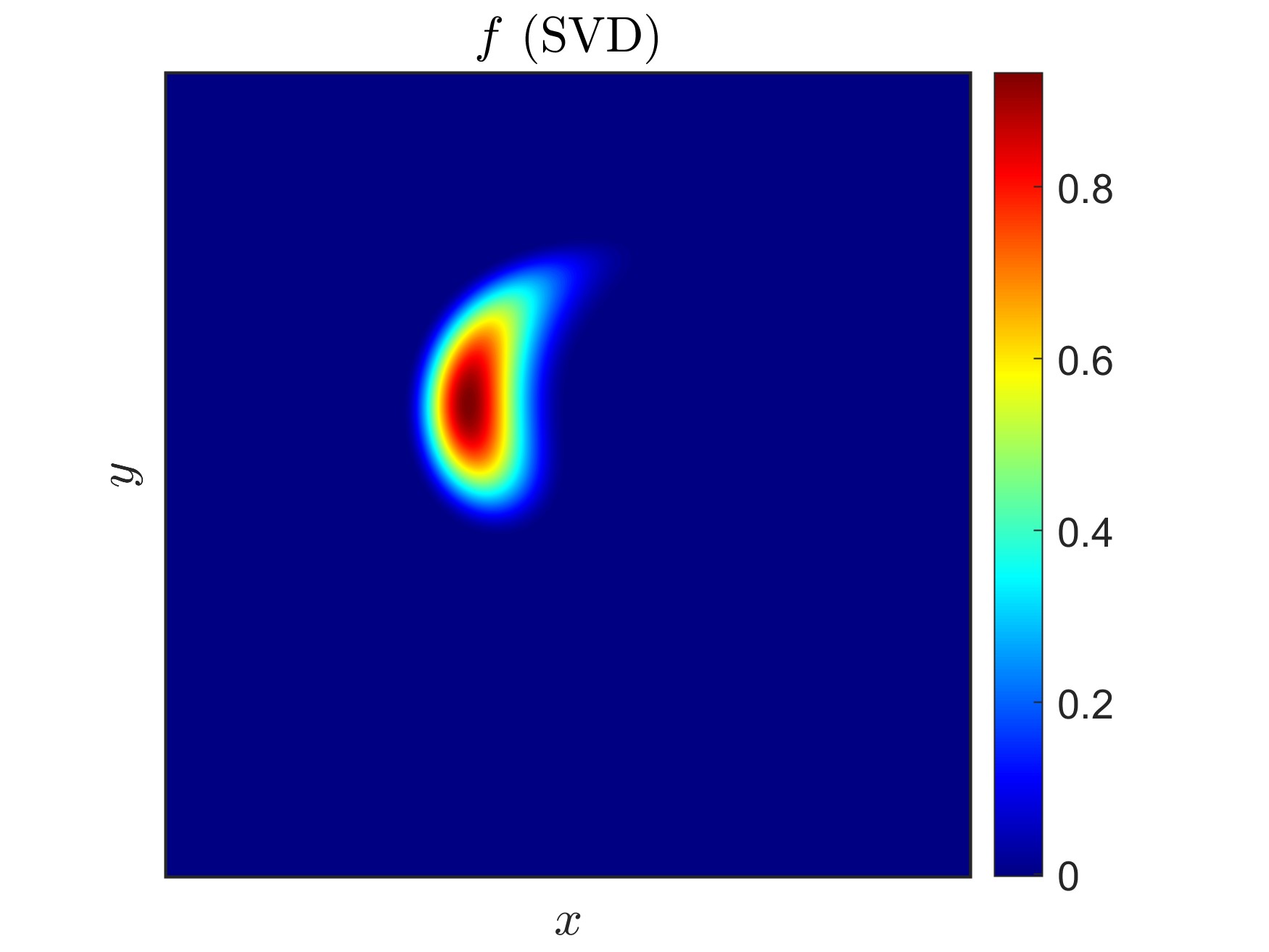}
  }

  \subfloat{
  \includegraphics[width=0.4\textwidth]{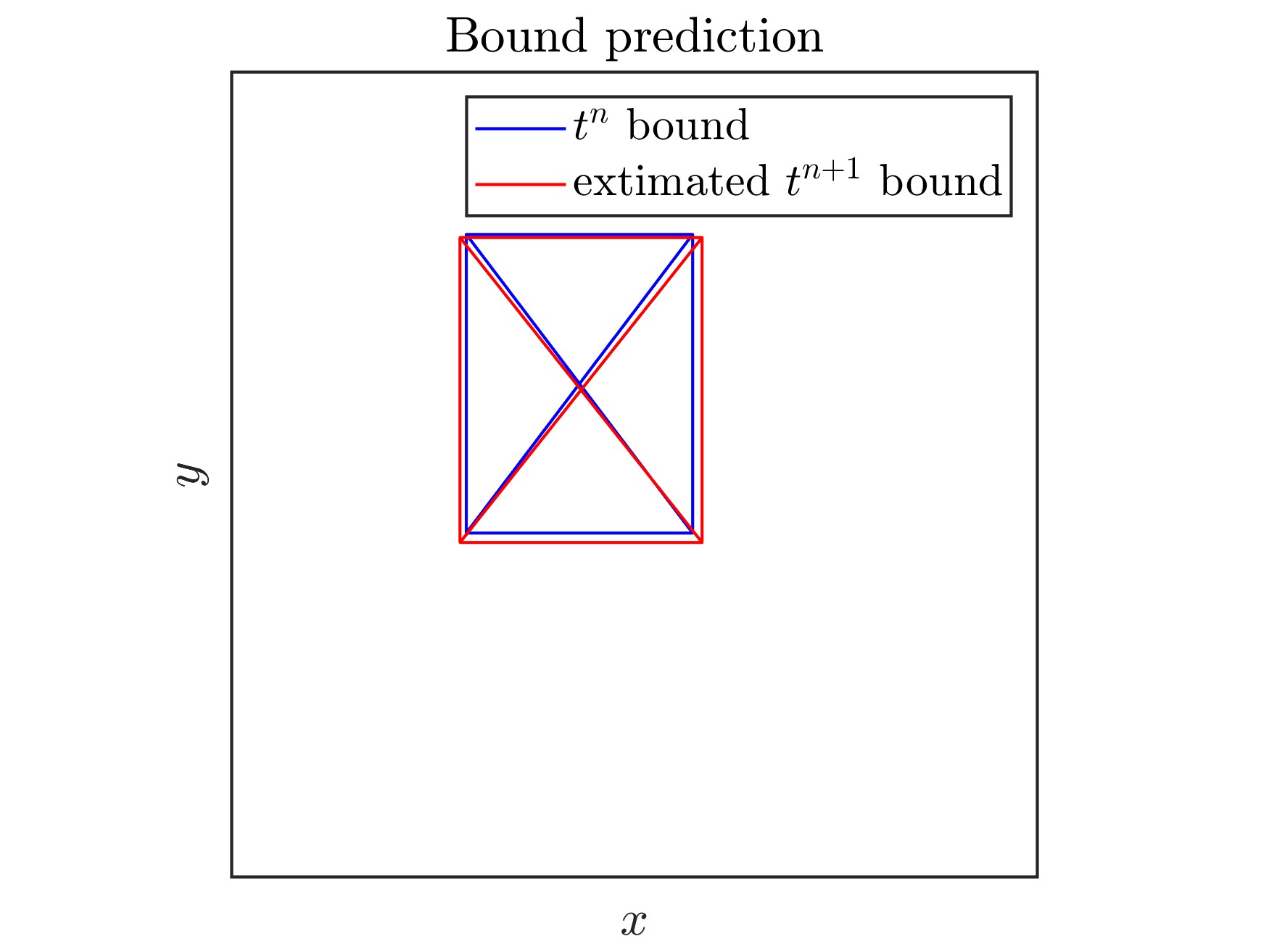}
  }
  \subfloat{
  \includegraphics[width=0.4\textwidth]{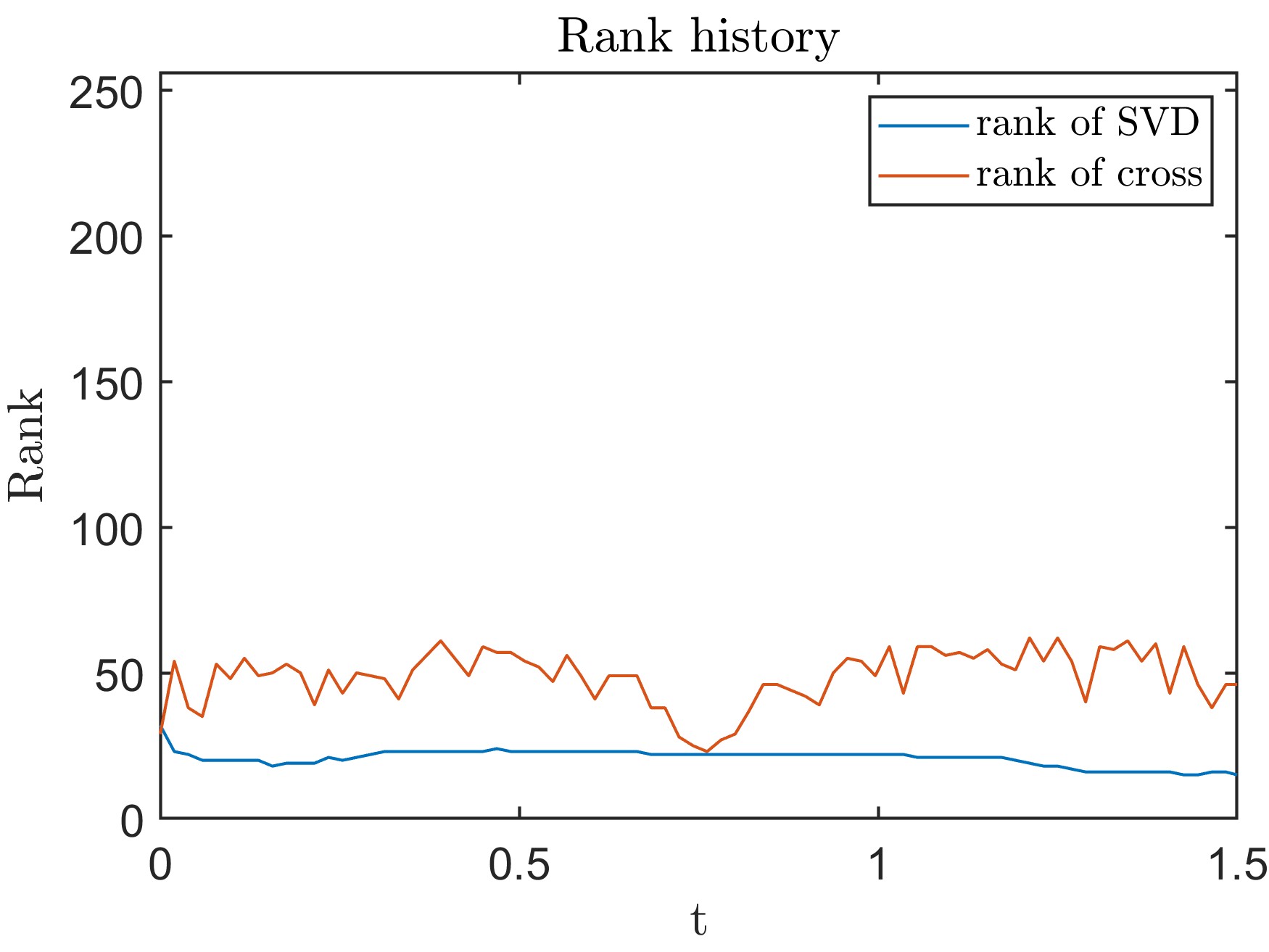}
  }
 \caption{(Swirling deformation flow). Top left: selected columns and rows for cross approximation at the half period. Top right: contour plot of the SVD solution at half period. Bottom left: visualization of the bound prediction procedure at the half period. Bottom right: rank history over the full period. A mesh of $256\times256$ and a CFL of 10 are used.}
  \label{fig:SDF}
\end{figure}
\end{exa}

\subsection{Nonlinear VP system}\label{sec:numer_VP}
In this subsection, we present two benchmark tests for the 1D1V VP system: the Landau damping and the bump-on-tail instability problems. Through these tests, we aim to verify the properties we examined for the linear advection equations, as well as further investigate the effectiveness of the RK exponential integrator \eqref{eq:CF3} and the enforcement of mass conservation. For all the tests in the subsection, we apply periodic boundary condition in the x-dimension and zero-boundary condition in the v-dimension.

\begin{exa}(Landau damping). Consider the 1D1V VP system with the initial condition
\begin{equation}\label{landau_damping}
f(x,v,t=0) = \frac1{\sqrt{2\pi}}\left(1+\alpha\cos(kx)\right)\exp\left(-\frac{v^2}{2}\right),\quad x\in[0,4\pi],\quad v\in[-2\pi,2\pi],
\end{equation}
where \(k = 0.01\), \(\alpha = 0.01\) for the weak Landau damping and \(k=0.5\), \(\alpha=0.5\) for the strong Landau damping.

A standard test for Landau damping is to verify the exponential decay or growth rate of the electric field. For weak Landau damping, a theoretical decay rate of -0.1533 is known. We use a \(256\times 256\) mesh and a CFL number of 10 to evaluate the discrete \(L^2\) norm of the electric field.
On the left side of \Cref{fig:electric_history_WLD_SLD}, the semi-log plot of the time evolution of the \(L^2\) norm of the electric field for weak Landau damping reflects the correct decay rate. For strong Landau damping, a similar result is shown on the right side of \Cref{fig:electric_history_WLD_SLD}. We compute the initial decay rate using the first two peaks of data and the growth rate using the ninth and eleventh peaks. The resulting decay and growth rates are approximately -0.2910 and 0.0810, respectively, which are very close to the results reported in the literature \cite{cheng1976integration, rossmanith2011positivity}.

In \Cref{fig:contour_rank_history_SLD}, for the strong Landau damping, we present the contour plot of the SVD solution at \(t = 40\) (left), the selected columns and rows of the cross approximation at \(t = 40\) (middle), and the rank history of the simulation from \(t = 0\) to \(t = 40\) (right). A mesh of \(256\times 256\) and a CFL number of 10 are used in this simulation. As shown, the SVD solution with a rank of around 40 effectively captures the filamentation structure of the strong Landau damping. The selected columns and rows exhibit a symmetrical structure, which aligns with the property of the real solution. Regarding the rank history, we observe a low-rank behavior throughout this limited-time simulation.

In \Cref{fig:mass_momentum_energy_SLD}, we present the time evolution of the relative deviation or deviation of discrete mass, momentum, and total energy. We adjust the computational range for $v$ to $[-10,10]$ for this simulation to prevent truncation error exceeding machine precision at the $v-$boundary. As shown, mass conservation is achieved, and the magnitudes of momentum and energy deviations are reasonable given our tolerance settings, i.e., \(\varepsilon_C = 10^{-4}\) and \(\varepsilon_S = 10^{-3}\).

Similar to the rigid body rotation problem, we investigate the temporal order of accuracy by fixing the spatial mesh and varying the CFL number. For strong Landau damping, we use two fixed meshes (\(128\times 128\) and \(256\times 256\)) and 20 different CFL numbers ranging from 1 to 100. The final simulation time is set to \(T = 5\). The reference solutions are computed using the fourth-order SL finite volume method \cite{zheng2022fourth} with a dense mesh of \(512\times 512\) and a CFL number of 1. We observe third-order accuracy in time for both the \(L^1\) and \(L^{\infty}\) errors, as shown in \Cref{fig:CFL_SLD}. Note that the temporal order accounts for errors from both the third-order RK exponential integrator \eqref{eq:CF3} and the S-stable third-order DIRK method \eqref{eq:S-stable_DIRK3}. 

For the strong Landau damping, in \Cref{fig:aver_rank_computing_time_SLD}, we display the rank behavior (left) and compression ratio of DOFs (right), similar to those for the rigid body rotation problem. In the left plot of \Cref{fig:iteration_number_SLD}, we also display the CFL number versus the average GMRES iteration count for the implicit solver \eqref{eq:S-stable_DIRK3}, both with and without the incomplete LU preconditioner, using a tolerance of \(10^{-14}\).
 The incomplete LU preconditioner significantly reduces the iteration numbers across different CFL numbers. The right plot of \Cref{fig:iteration_number_SLD} presents the computing time versus the grid point per dimension $N$ with a fixed time step size \(\Delta t = 0.01\) and final time \(t = 5\). We observe a linear complexity with respect to $N$, due to a complexity analysis of \(\mathcal{O}(Nr)\) flops for the local SL-FD solvers, \(\mathcal{O}(Nr^2 + r^3)\) flops for the SVD truncation, and \(\mathcal{O}(N)\) flops for the implicit update of the charge density.

\begin{figure}[!htbp]
\centering
  \subfloat{
  \includegraphics[width=0.4\textwidth]{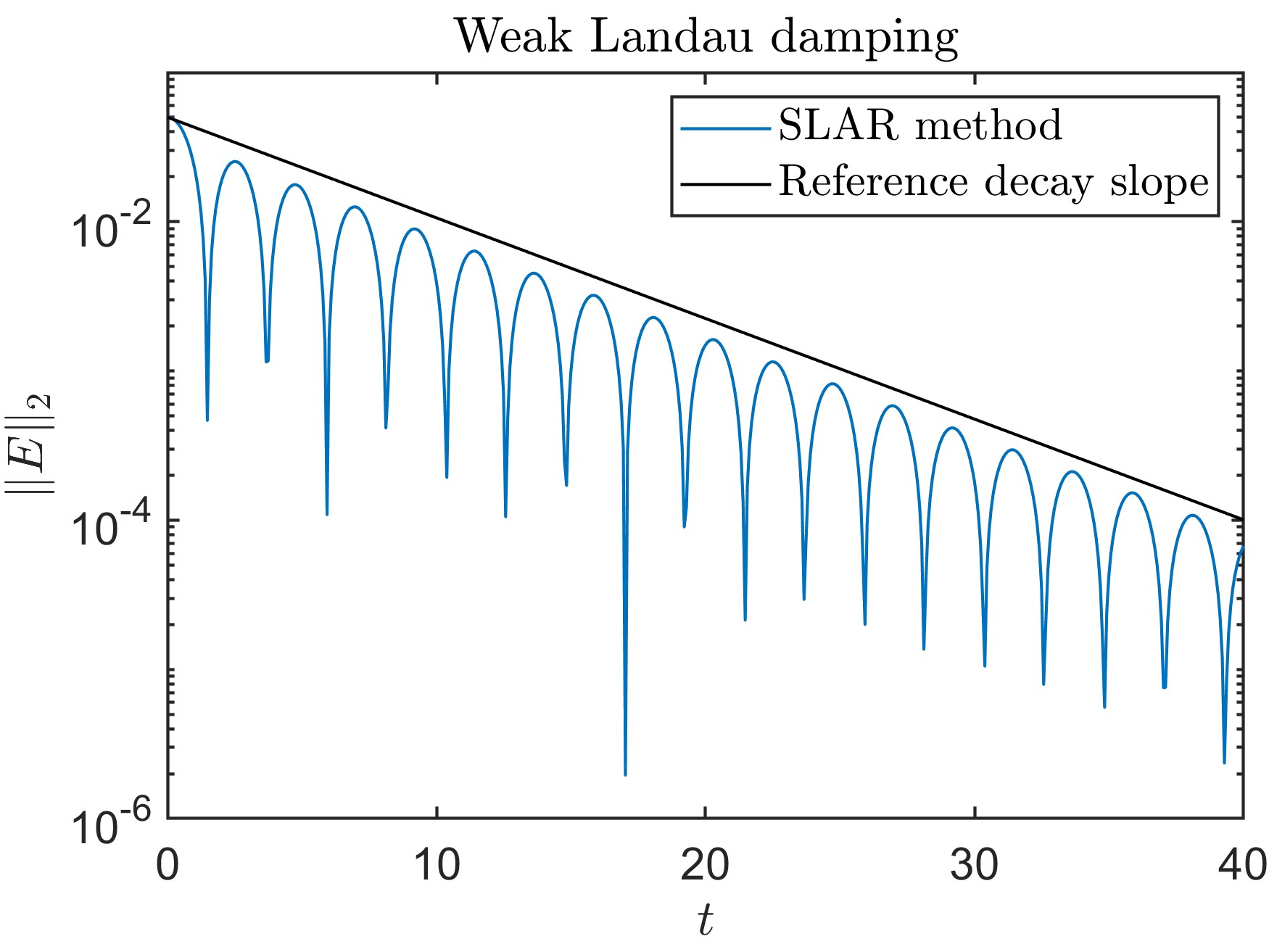}
  }
  \subfloat{
  \includegraphics[width=0.4\textwidth]{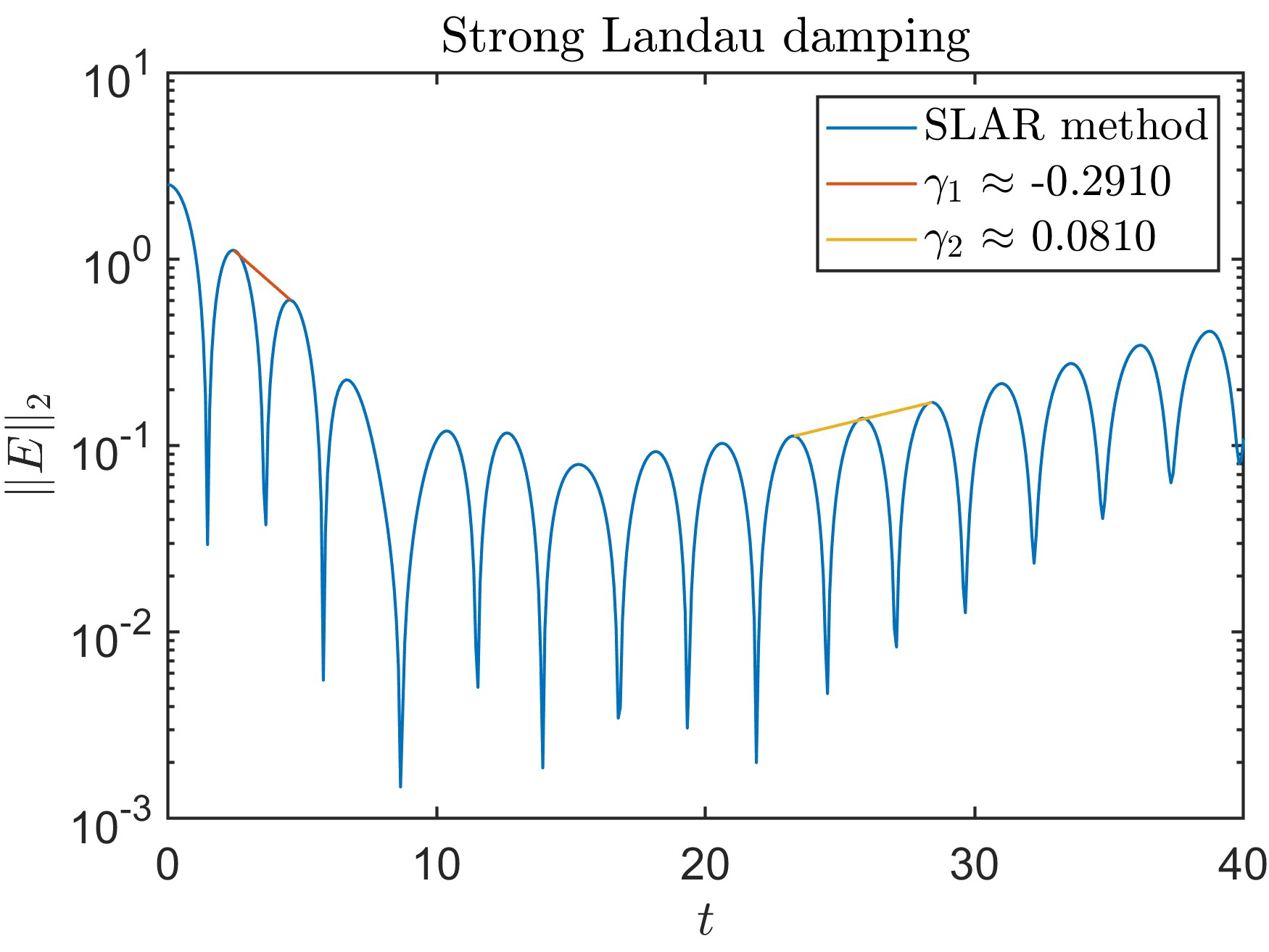}
  }
 \caption{(Landau damping). Time evolution of the $L^2$ norm of the electric filed for the weak (left) and strong (right) Landau dampings with a mesh of $256\times256$ and a CFL of 10. For the weak Landau damping, $\varepsilon_C = 10^{-5}$ and $\varepsilon_S = 10^{-4}$ are used to match the resolution requirement.}
  \label{fig:electric_history_WLD_SLD}
\end{figure}

\begin{figure}[!htbp]
\centering
  \subfloat{
  \includegraphics[width=0.32\textwidth]{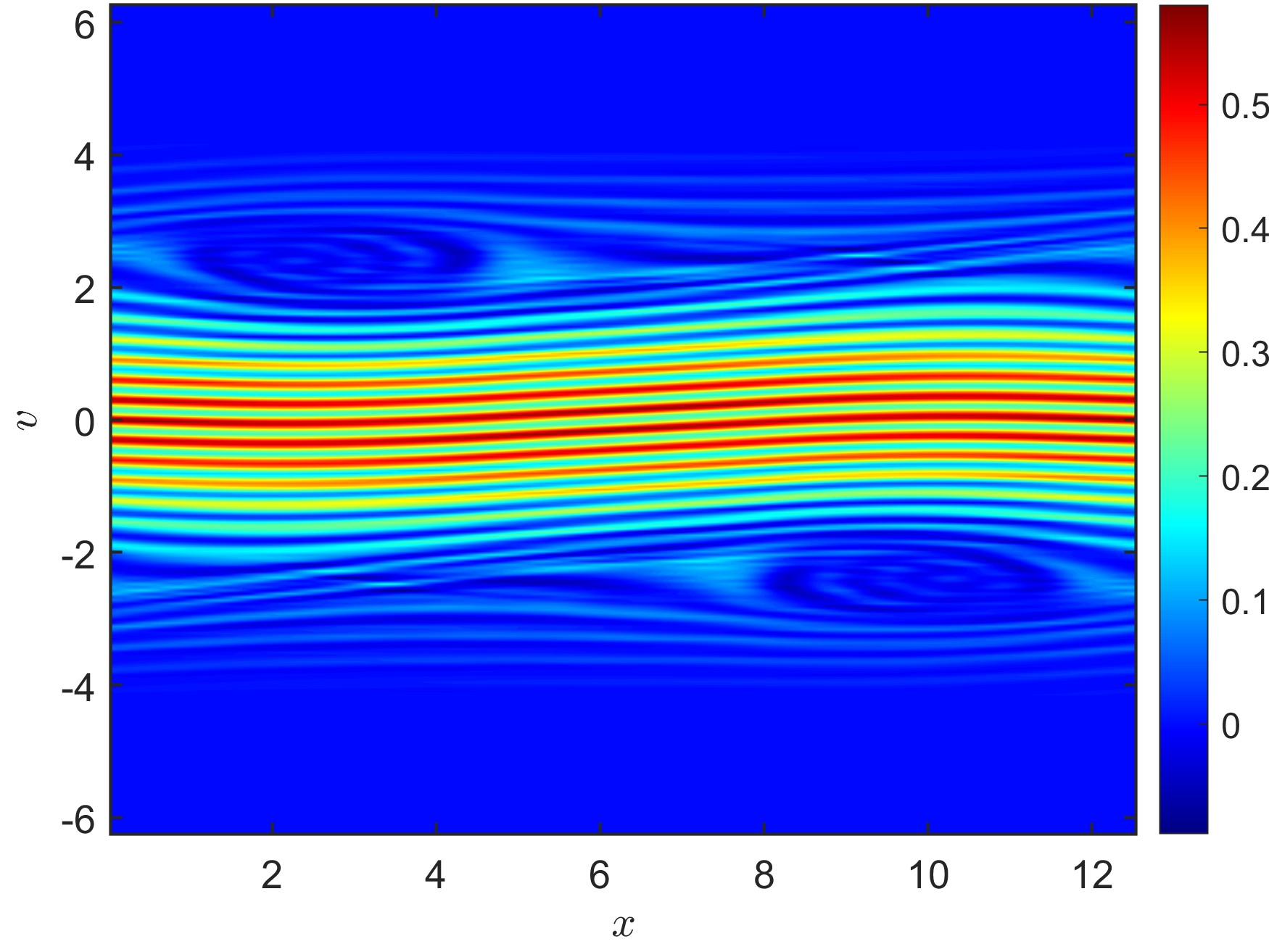}
  }
  \subfloat{
  \includegraphics[width=0.32\textwidth]{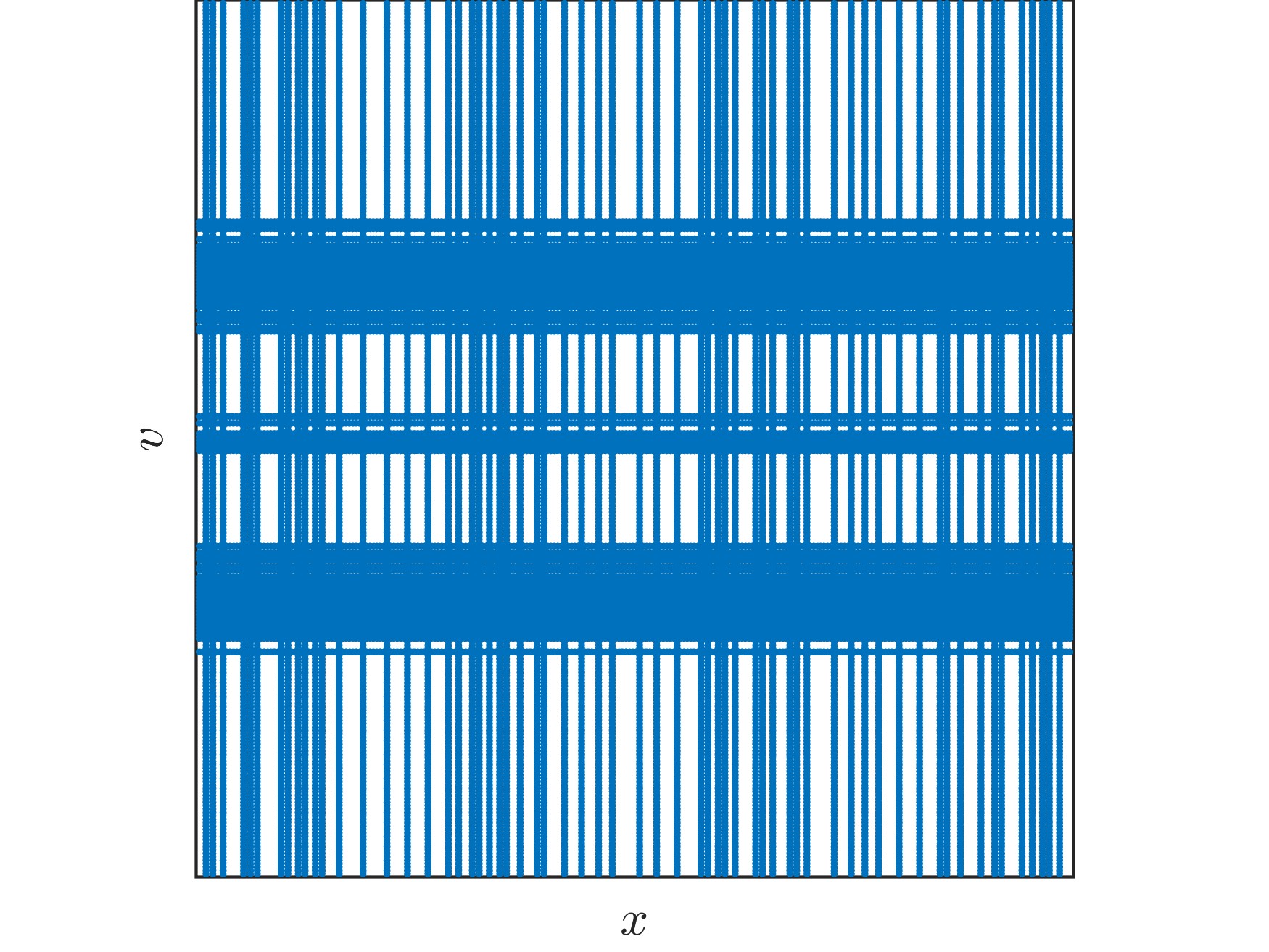}
  }
  \subfloat{
  \includegraphics[width=0.32\textwidth]{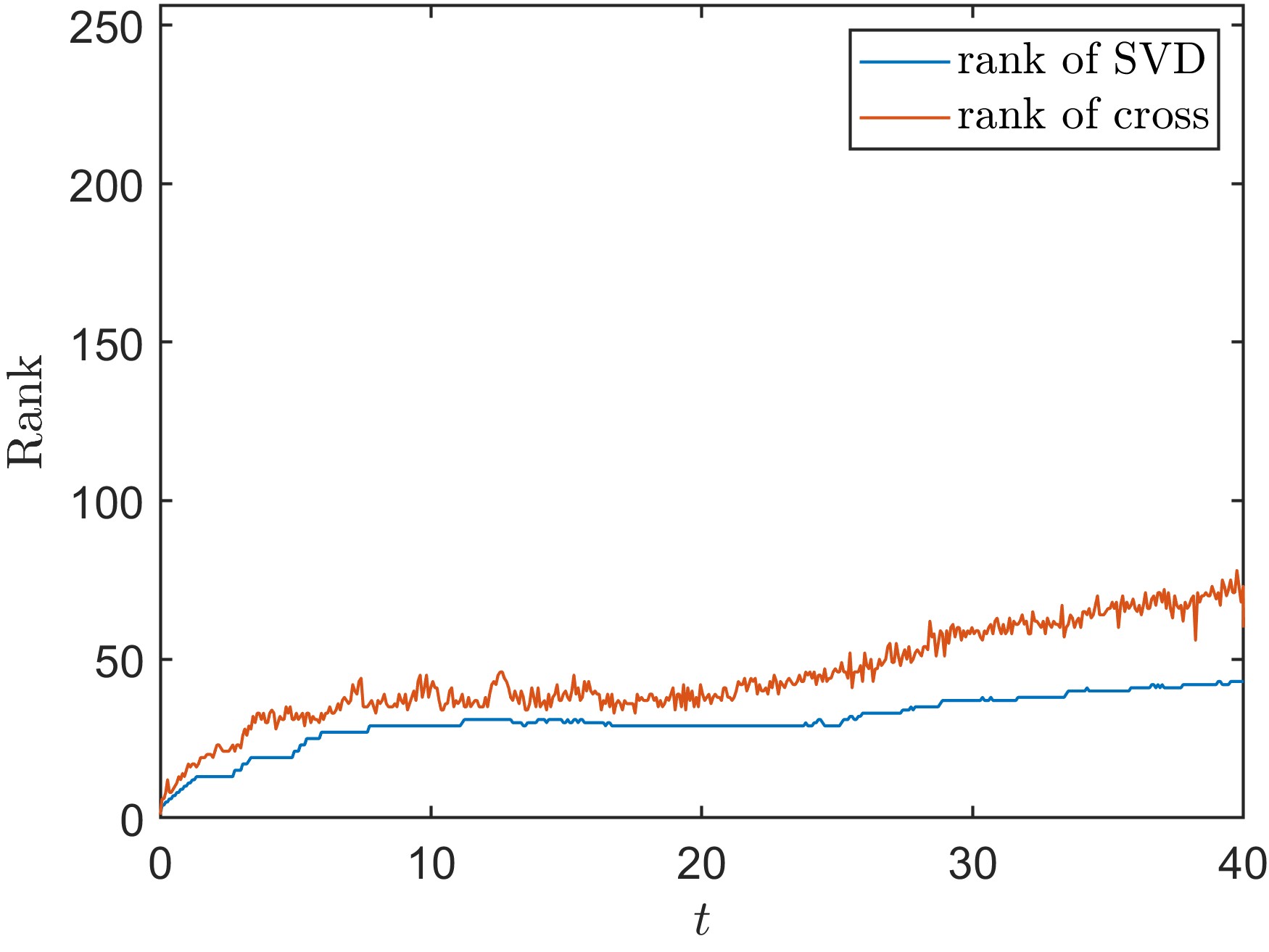}
  }
 \caption{(Strong Landau damping). Left: contour plot of the numerical solution at $t = 40$. Middle: the selected columns and rows at $t = 40$. Right: rank history of the simulation from $t=0$ to $40$. The simulation uses a mesh of $256\times256$ and a CFL of 10.}
  \label{fig:contour_rank_history_SLD}
\end{figure}

\begin{figure}[!htbp]
\centering
  \subfloat{
  \includegraphics[width=0.32\textwidth]{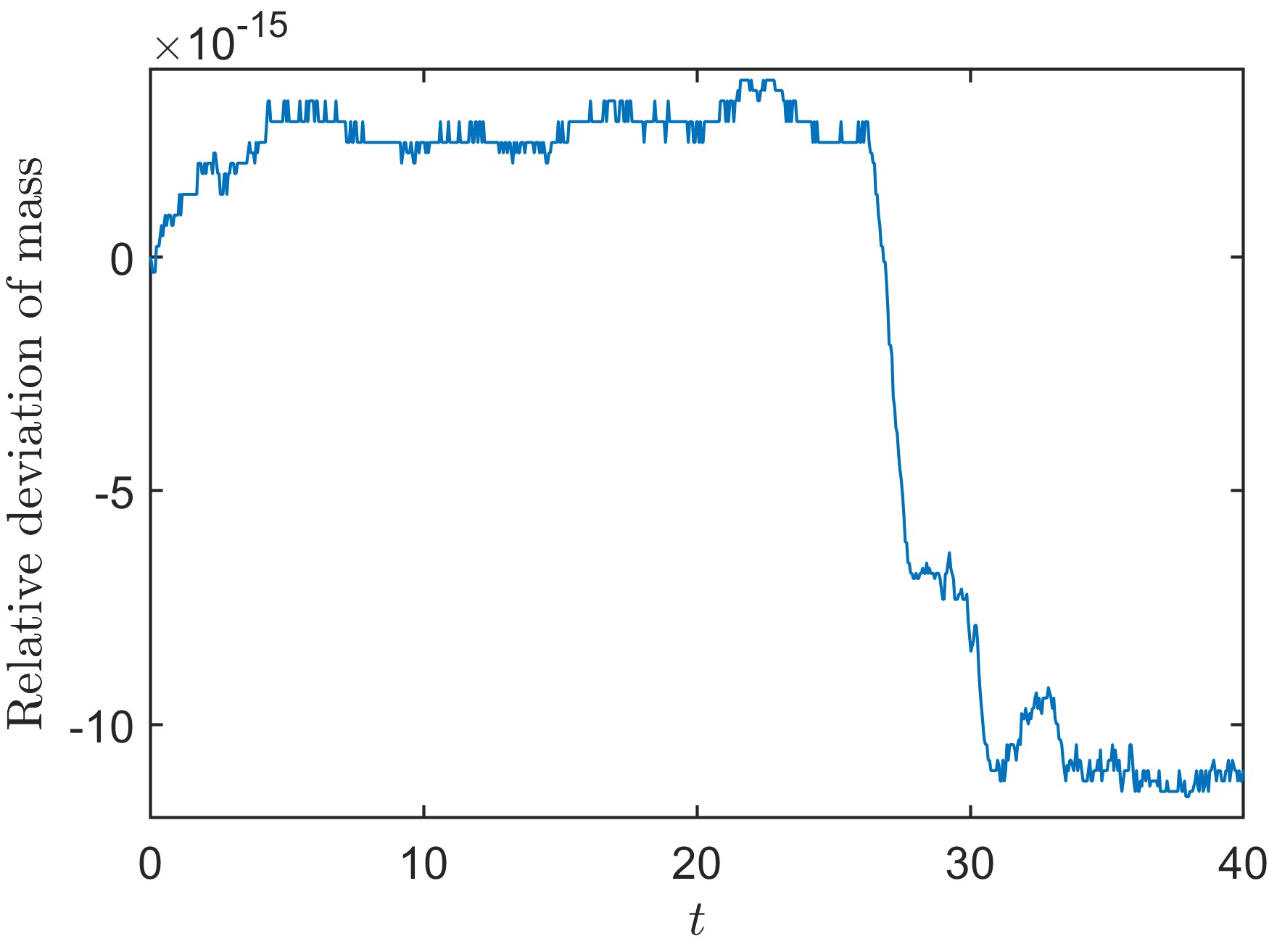}
  }
  \subfloat{
  \includegraphics[width=0.32\textwidth]{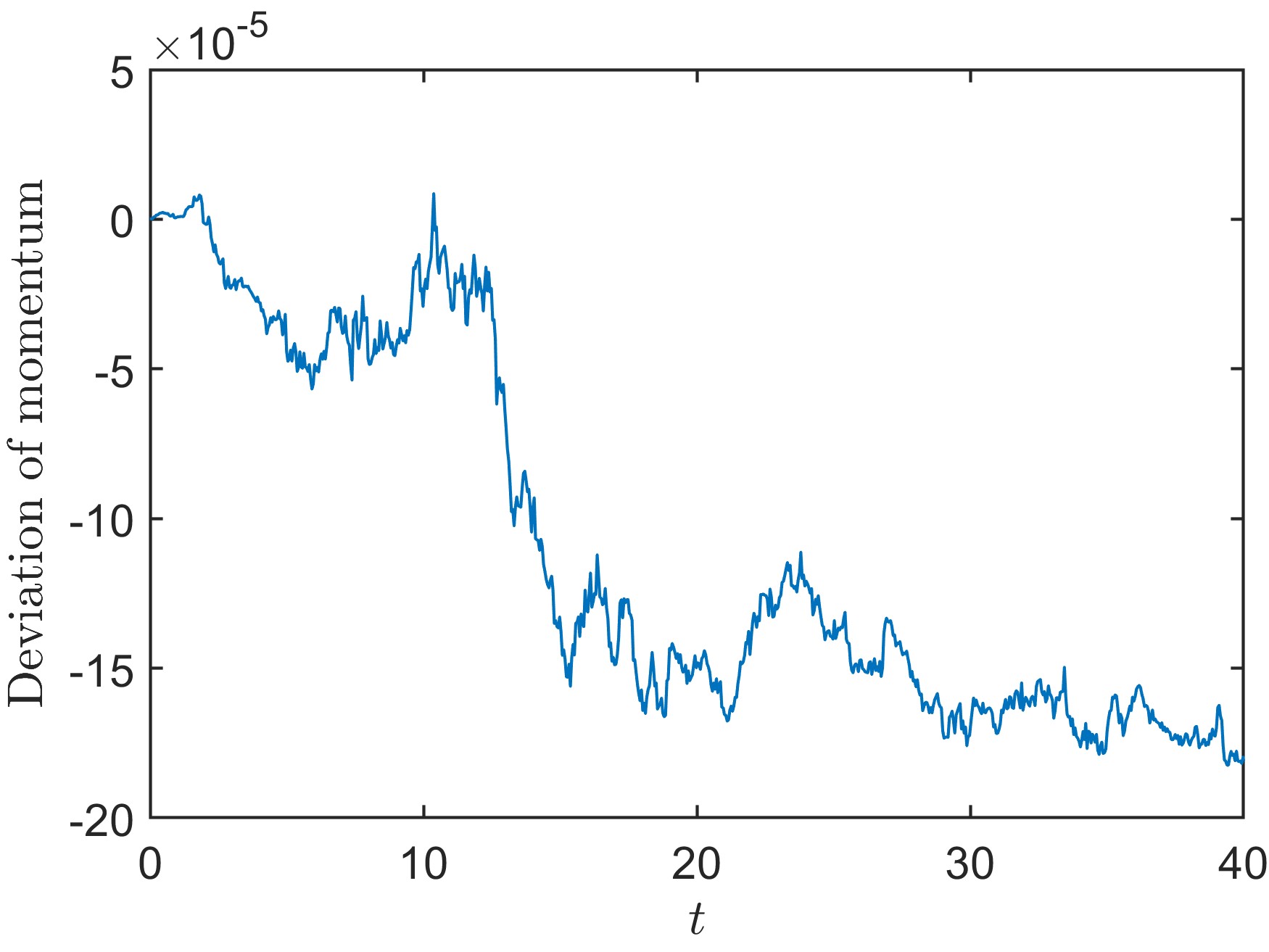}
  }
  \subfloat{
  \includegraphics[width=0.32\textwidth]{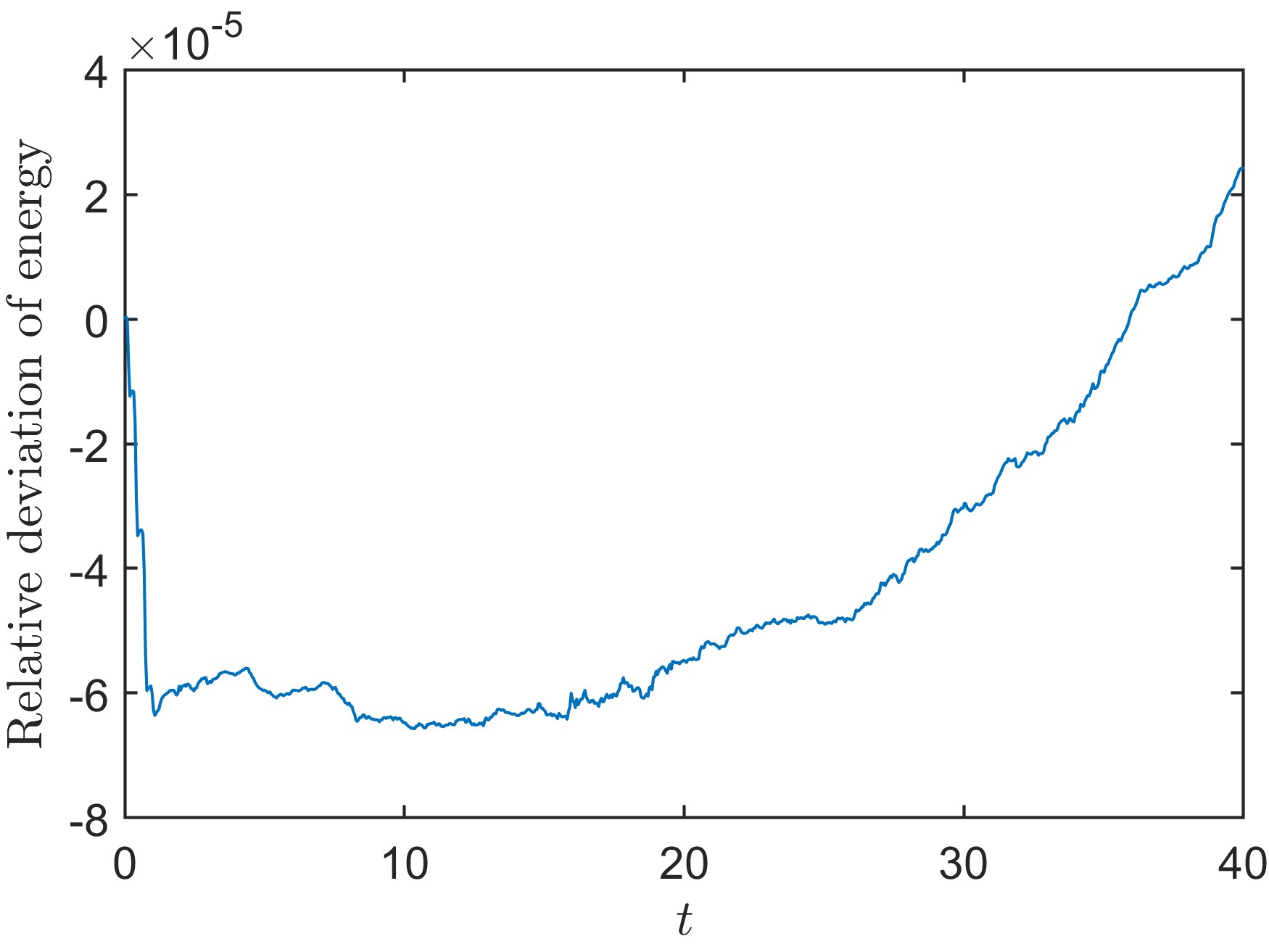}
  }

 \caption{(Strong Landau damping). performance of preserving mass, momemtum and energy of the simulation with a mesh of $256\times256$ and a CFL of 10. The computational range for $v$ is adjusted to $[-10,10]$.}
  \label{fig:mass_momentum_energy_SLD}
\end{figure}

\begin{figure}[!htbp]
\centering
  \subfloat{
  \includegraphics[width=0.4\textwidth]{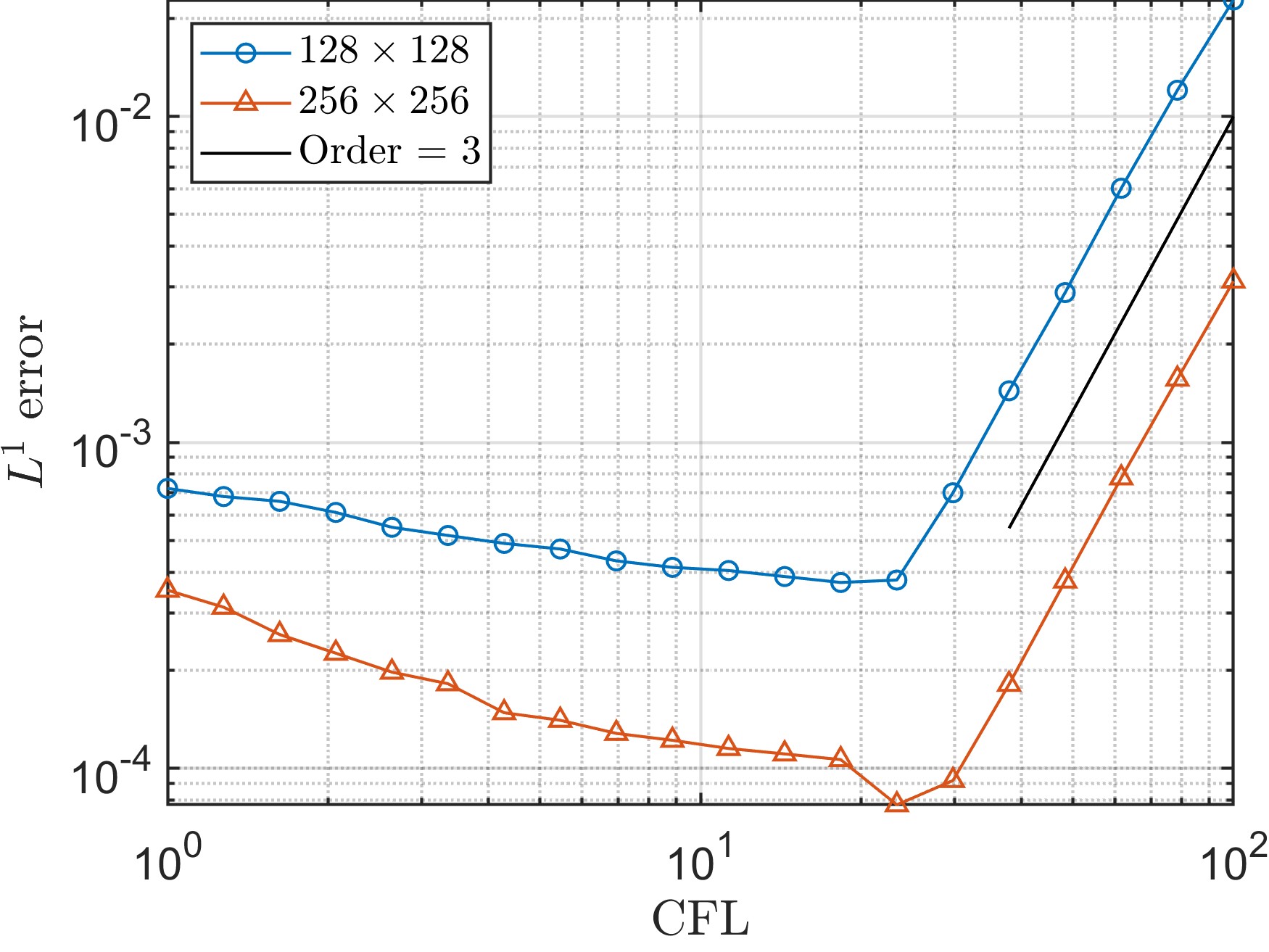}
  }
  \subfloat{
  \includegraphics[width=0.4\textwidth]{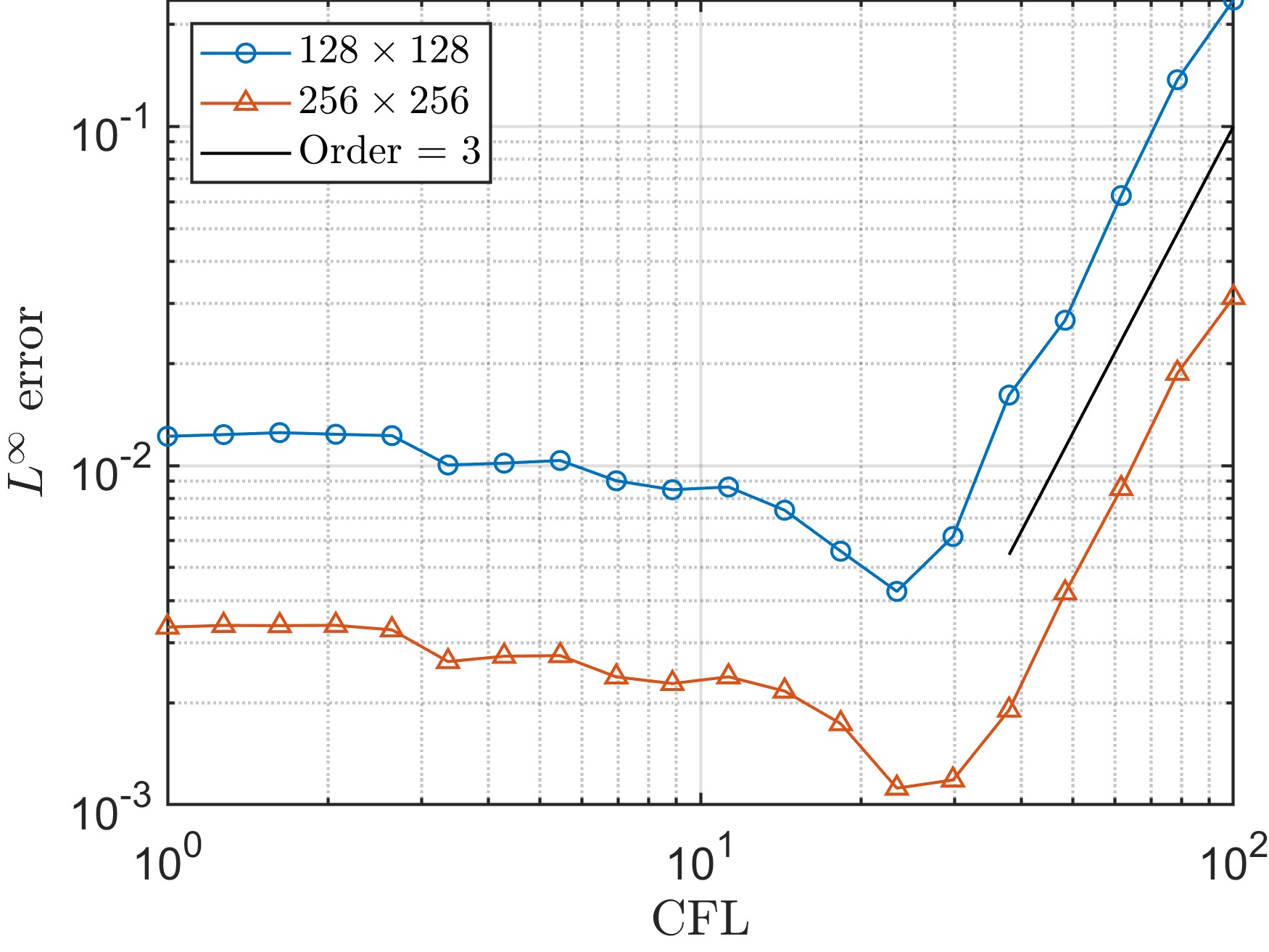}
  }
 \caption{(Strong Landau damping). Log-log plots of CFL numbers versus $L^1$ and $L^{\infty}$ errors with two sets of fixed meshes, $128\times128$ and $256\times256$ at $t = 5$. }
  \label{fig:CFL_SLD}
\end{figure}

\begin{figure}[!htbp]
\centering
  \subfloat{
  \includegraphics[width=0.4\textwidth]{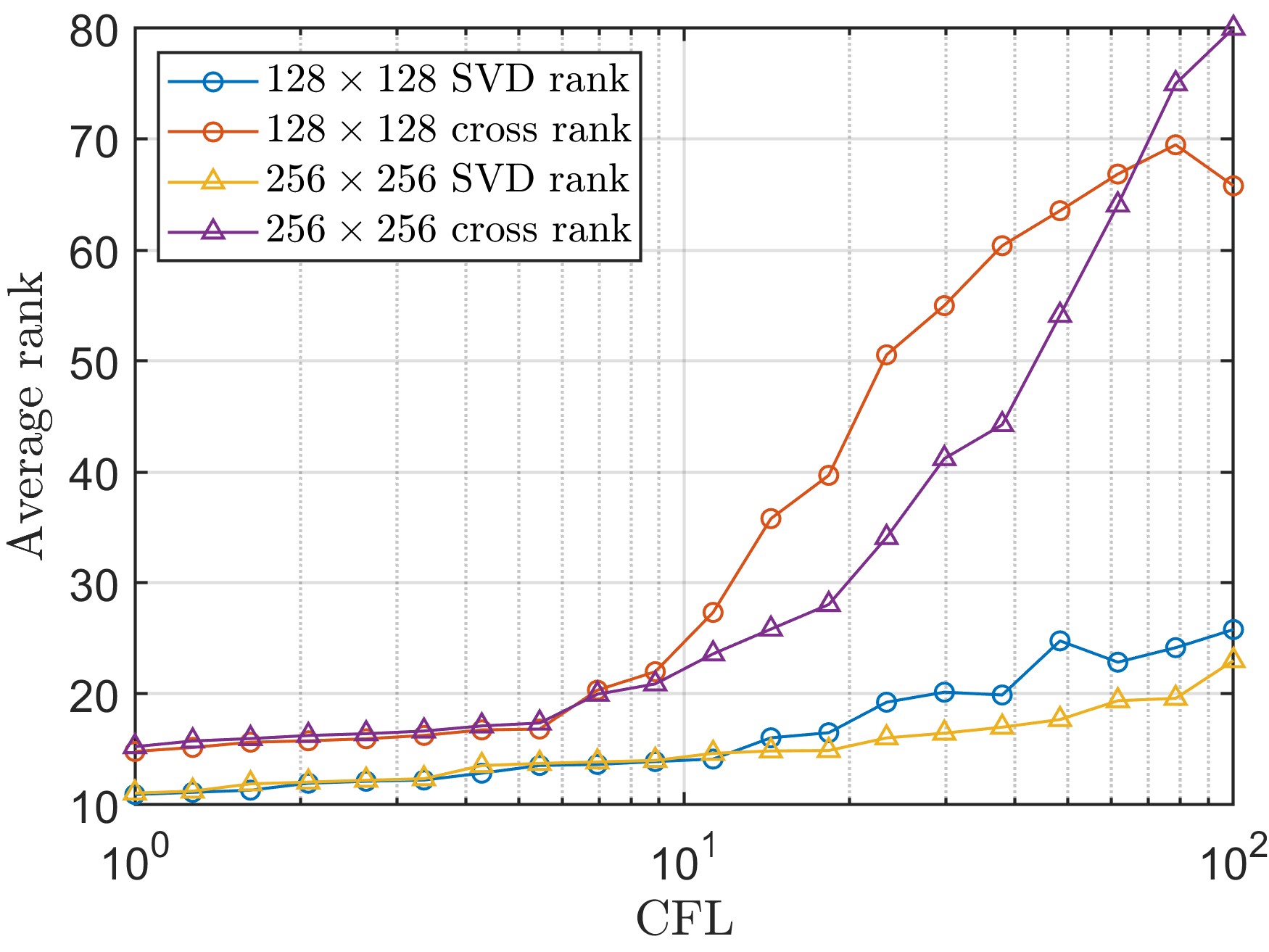}
  }
  \subfloat{
  \includegraphics[width=0.4\textwidth]{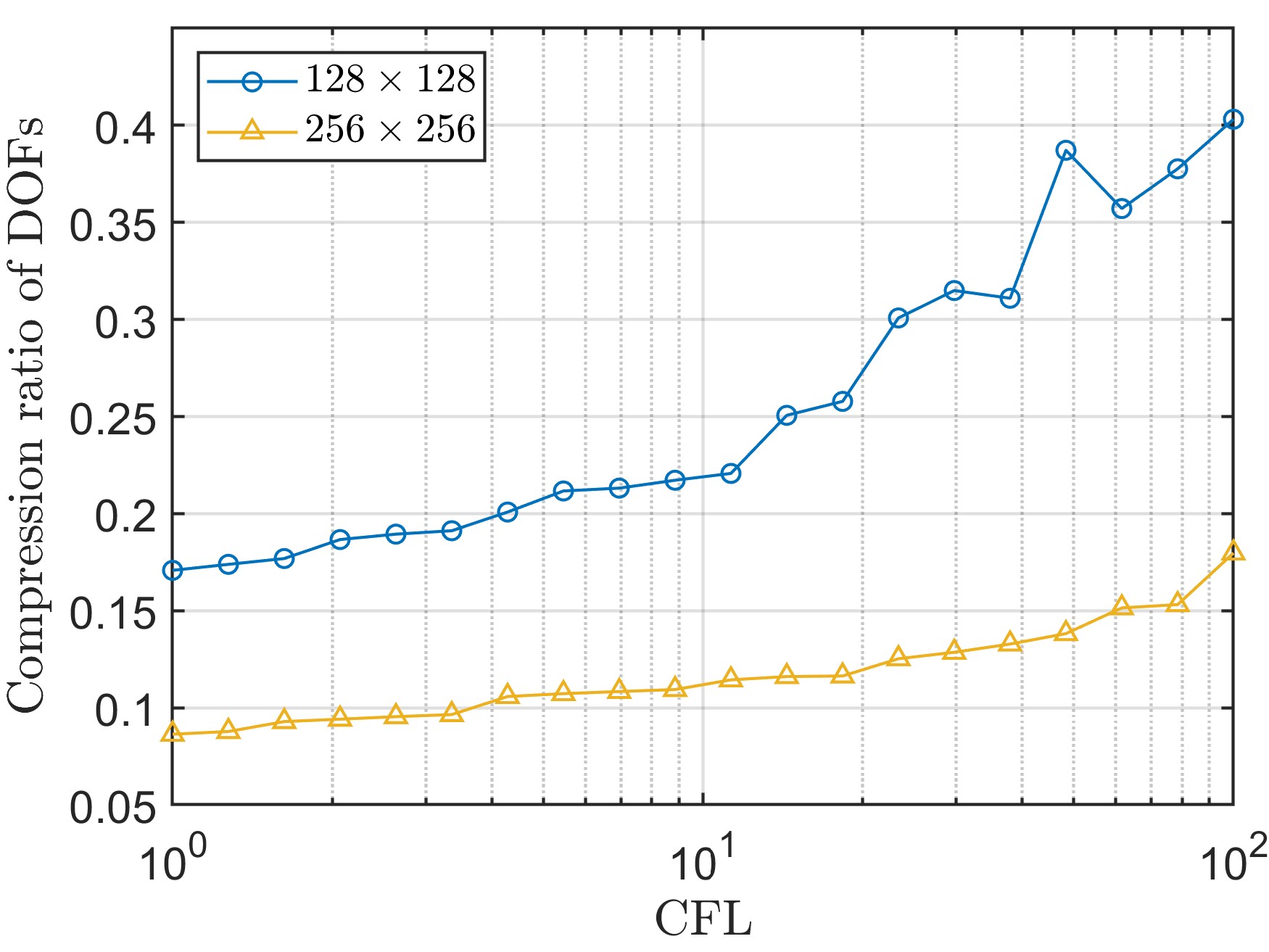}
  }
 \caption{(Strong Landau damping). Left: semi-log plot of CFL numbers versus average ranks of the simulations in \Cref{fig:CFL_SLD} (t = 5). Right: log-log plot of CFL numbers versus compression ratios of DOFs of the simulations in \Cref{fig:CFL_SLD} (t = 5).}
  \label{fig:aver_rank_computing_time_SLD}
\end{figure}

\begin{figure}[!htbp]
\centering
  \subfloat{
  \includegraphics[width=0.4\textwidth]{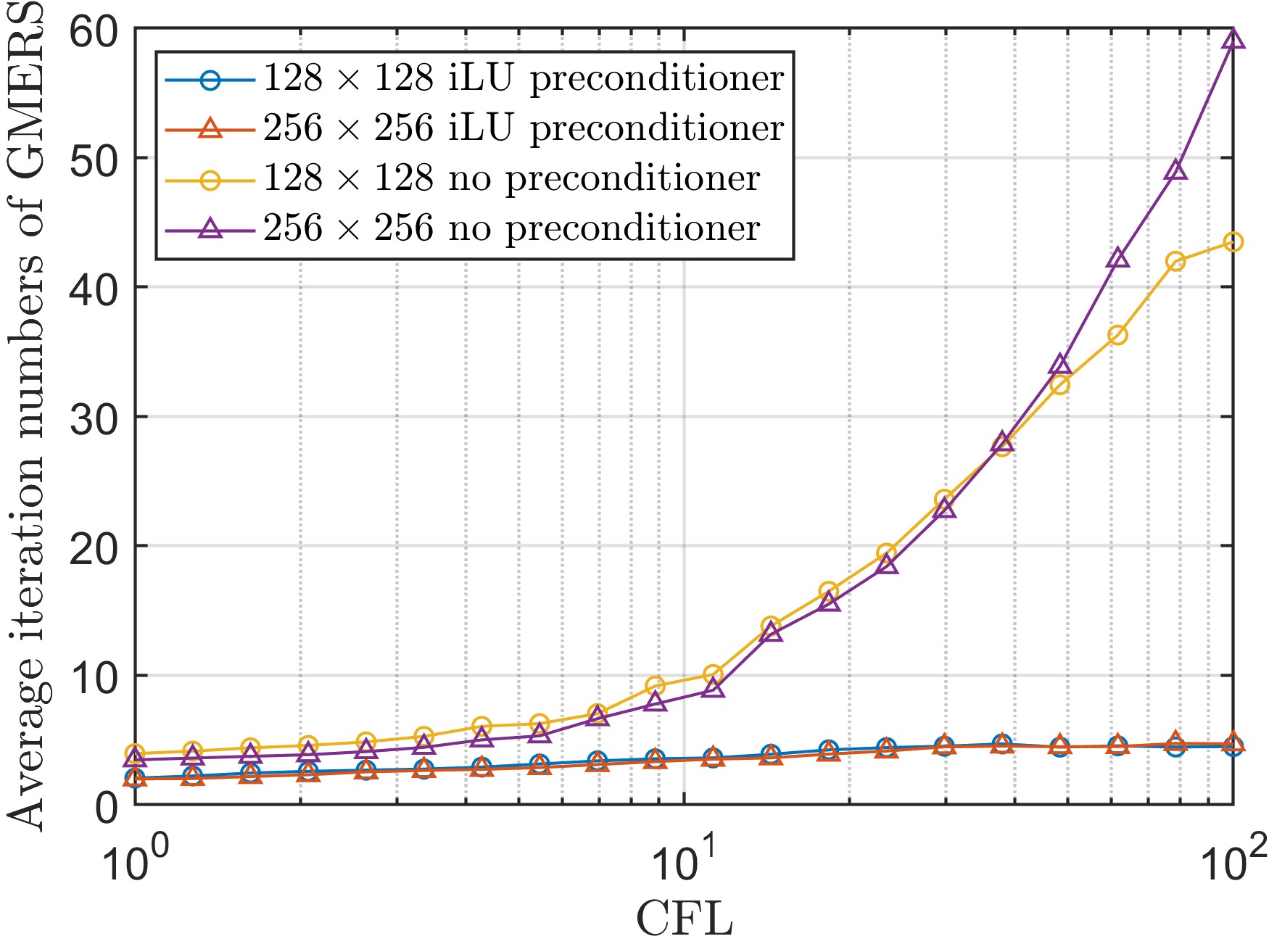}
  }
  \subfloat{
  \includegraphics[width=0.4\textwidth]{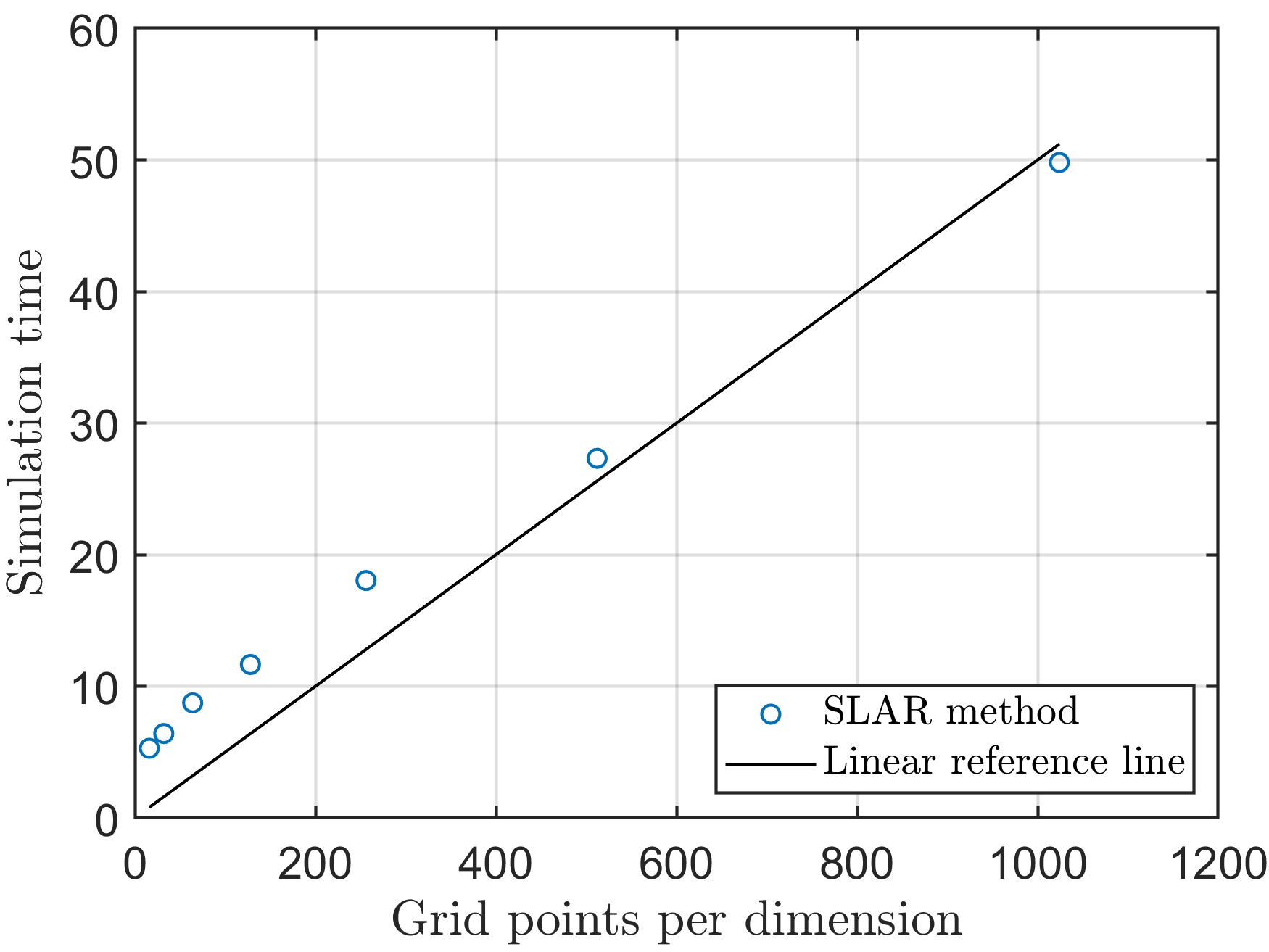}
  }
 \caption{(Strong Landau damping). Left: semi-log plot of CFL numbers versus average iteration numbers for the sparse GMRES method with a tolerance of $10^{-14}$ for the simulations in \Cref{fig:CFL_SLD} (t = 5). Right: plot of grid points per dimension versus simulation time with a fixed time step size $\Delta t = 0.01$ and for the target time $t = 5$. }
  \label{fig:iteration_number_SLD}
\end{figure}



\end{exa}

\begin{exa}(Bump-on-tail instability). Consider the bump-on-tail instability with the initial condition
\begin{equation}\label{eq:BOT}
\begin{split}
f(x,v,t=0) = (1+\alpha\cos(kx))\left(n_p\exp\left(-\frac{v^2}{2}\right) + n_b\exp\left(-\frac{(v-u)^2}{2v_t}\right)\right),\\
\quad x\in[0,20\pi/3],\quad v\in[-13,13],
\end{split}
\end{equation}
where \(\alpha = 0.04\), \(k = 0.3\), \(n_p=\frac{9}{10\sqrt{2\pi}}\), \(n_b=\frac{2}{10\sqrt{2\pi}}\), \(u=4.5\), and \(v_t=0.5\). In \Cref{fig:contour_rank_history_BOT}, we present the contour plot of the SVD solution at \(t = 40\) (left), the selected columns and rows of the cross approximation at \(t = 40\) (middle), and the rank history of the simulation from \(t = 0\) to \(t = 40\) (right) with  a \(256\times 256\) mesh and a CFL number of 10. The contour plot is consistent with existing results in the literature \cite{cai2018high}. The selected x-dimension slices are clustered near regions of rich information. Throughout this limited-time simulation, the rank history consistently exhibits low-rank behavior. On the top left of \Cref{fig:mass_momentum_energy_BOT}, we present the time evolution of the \(L^2\) norm of the electric field. The result closely matches those reported in the literature \cite{cai2018high,guo2024local}. On the top right and bottom of \Cref{fig:mass_momentum_energy_BOT}, we present the performance in preserving mass, momentum, and total energy; similar behavior is observed for the Landau damping case.

\begin{figure}[!htbp]
\centering
  \subfloat{
  \includegraphics[width=0.32\textwidth]{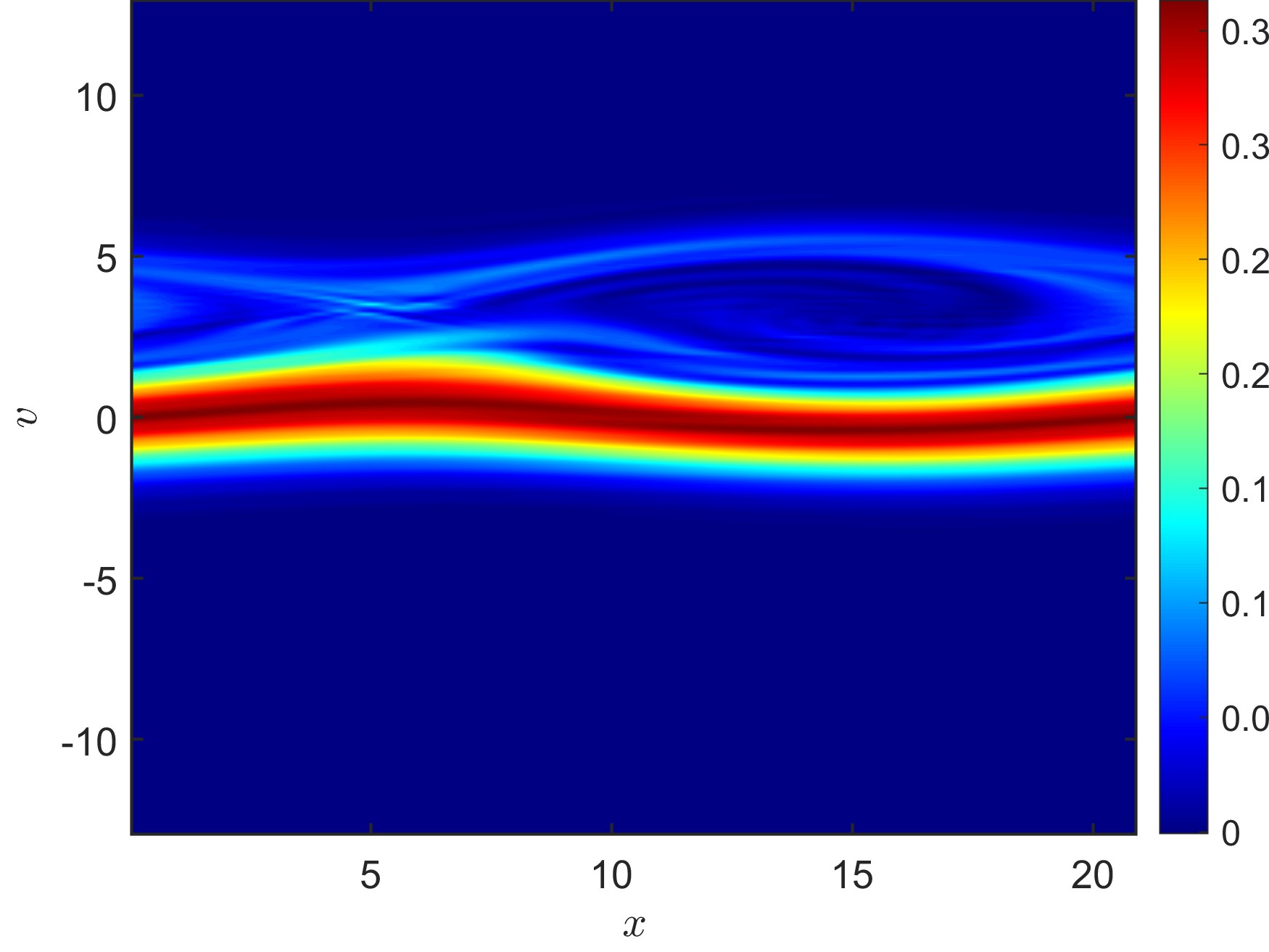}
  }
  \subfloat{
  \includegraphics[width=0.32\textwidth]{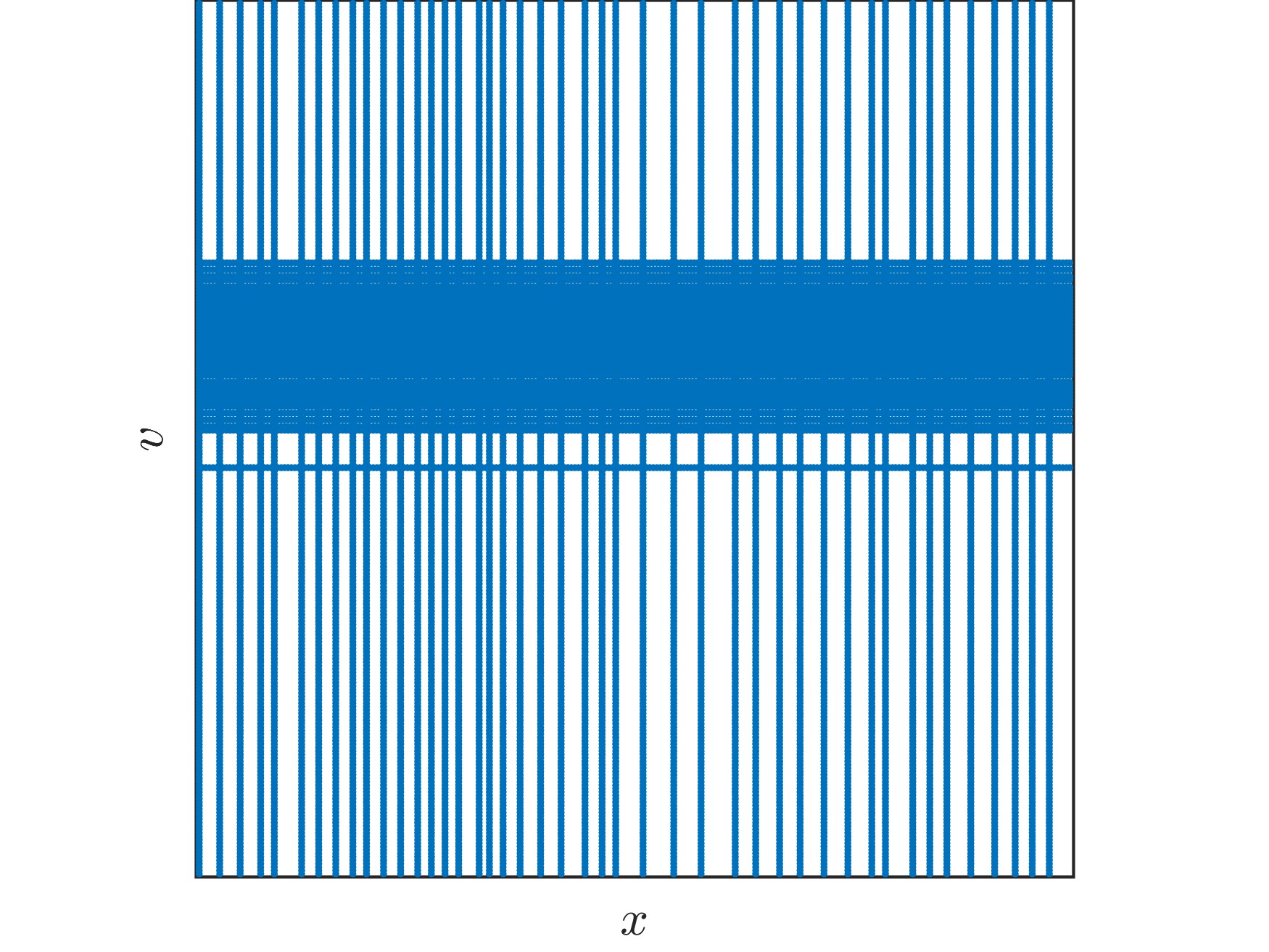}
  }
  \subfloat{
  \includegraphics[width=0.32\textwidth]{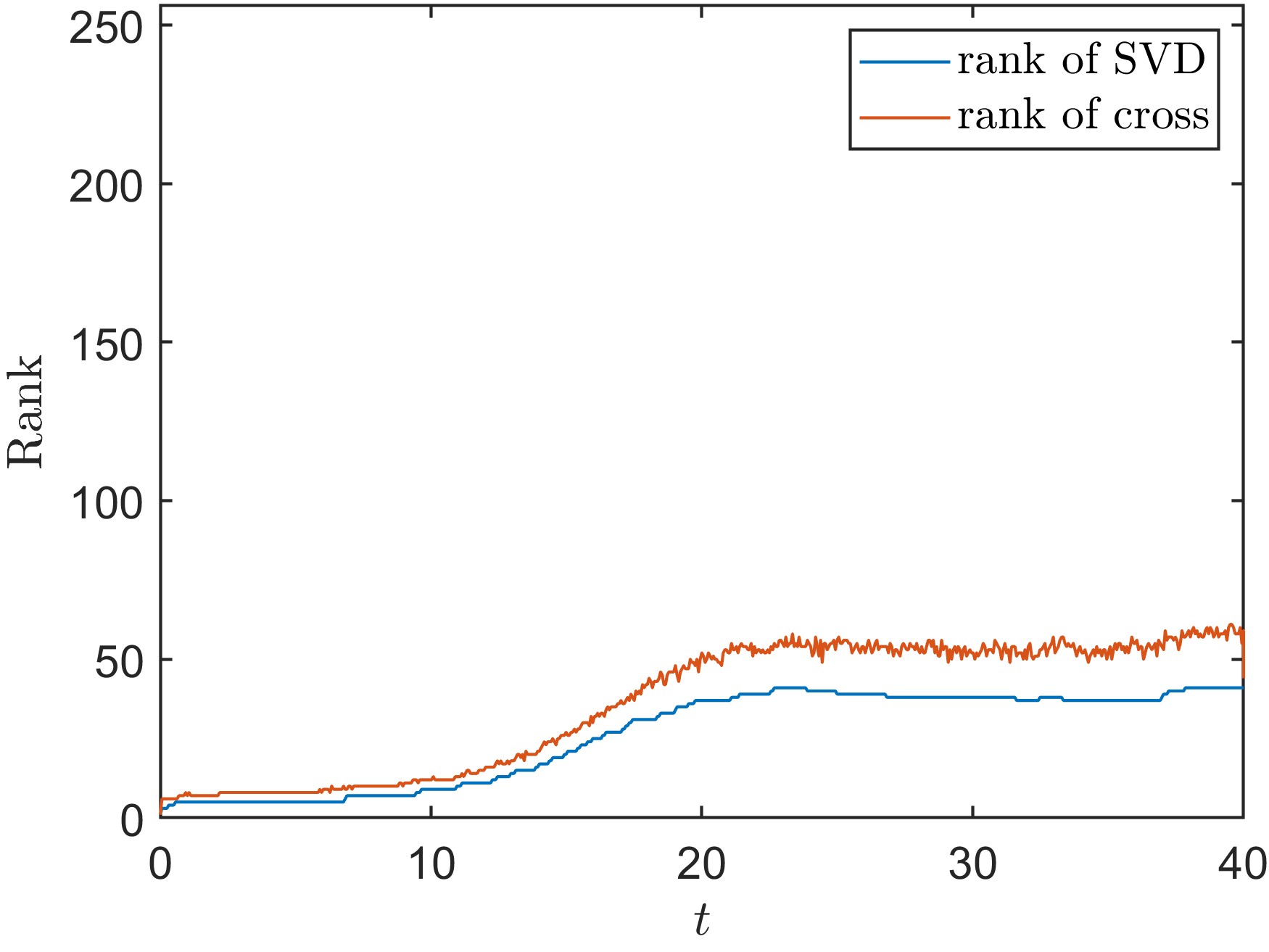}
  }
 \caption{(Bump-on-tail instability). Left: contour plot of the numerical solution at $t = 40$. Middle: the selected columns and rows at $t = 40$. Right: rank history of the simulation from $t=0$ to $40$. The simulation uses a mesh of $256\times256$ and a CFL of 10.}
  \label{fig:contour_rank_history_BOT}
\end{figure}

\begin{figure}[h]
\centering
  \subfloat{
  \includegraphics[width=0.36\textwidth]{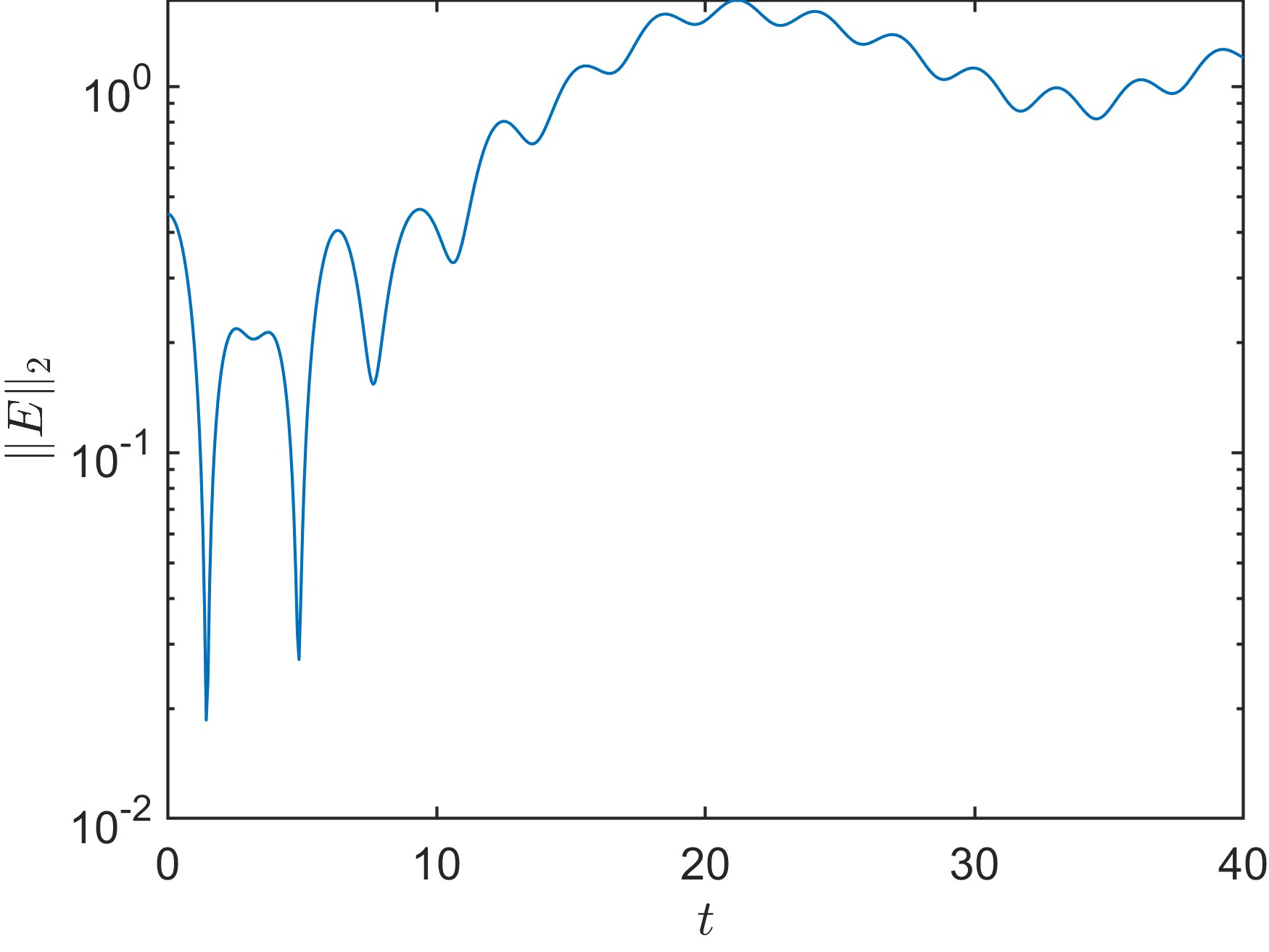}
  }
  \subfloat{
  \includegraphics[width=0.36\textwidth]{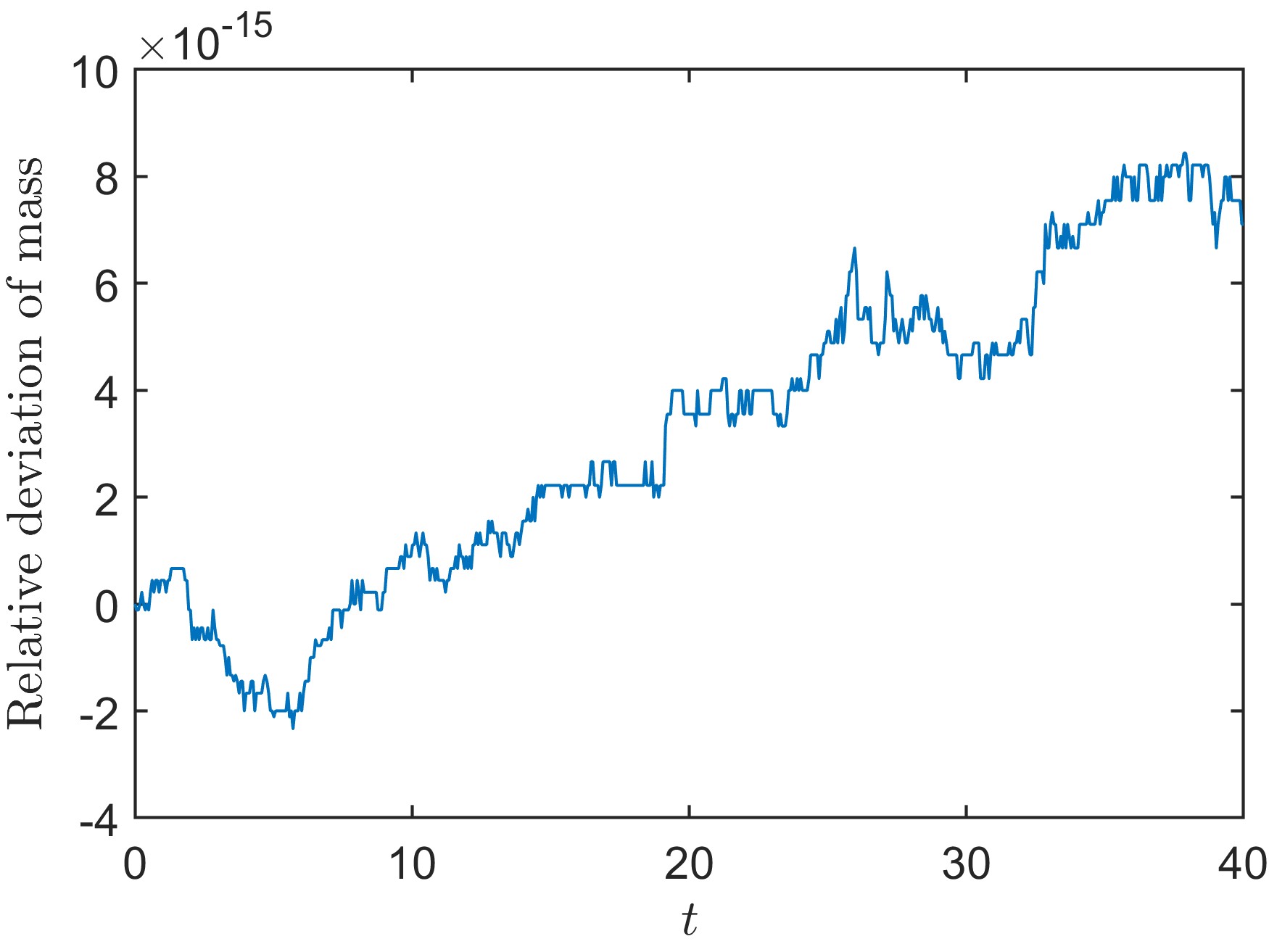}
  }
  
  \subfloat{
  \includegraphics[width=0.36\textwidth]{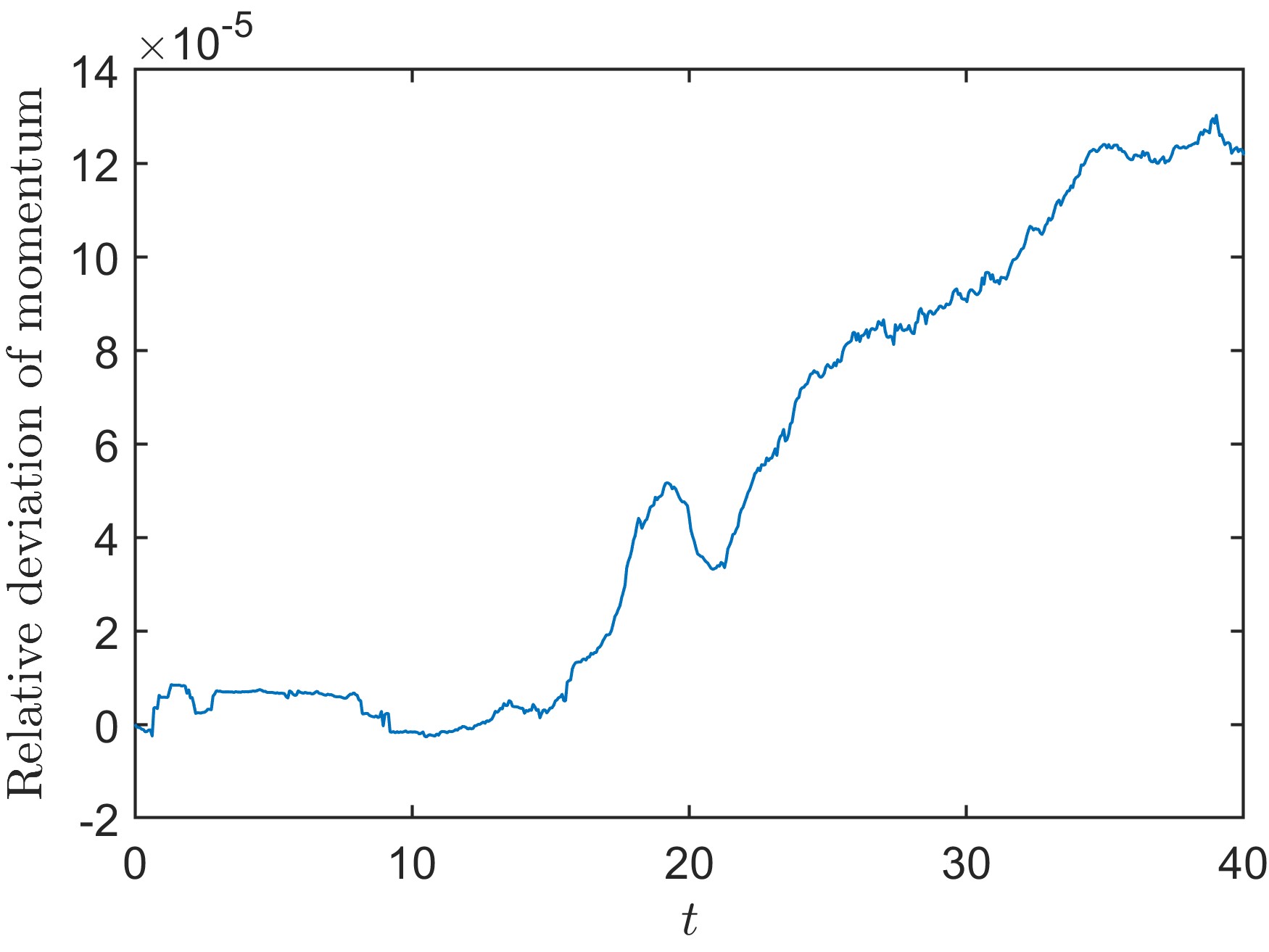}
  }
  \subfloat{
  \includegraphics[width=0.36\textwidth]{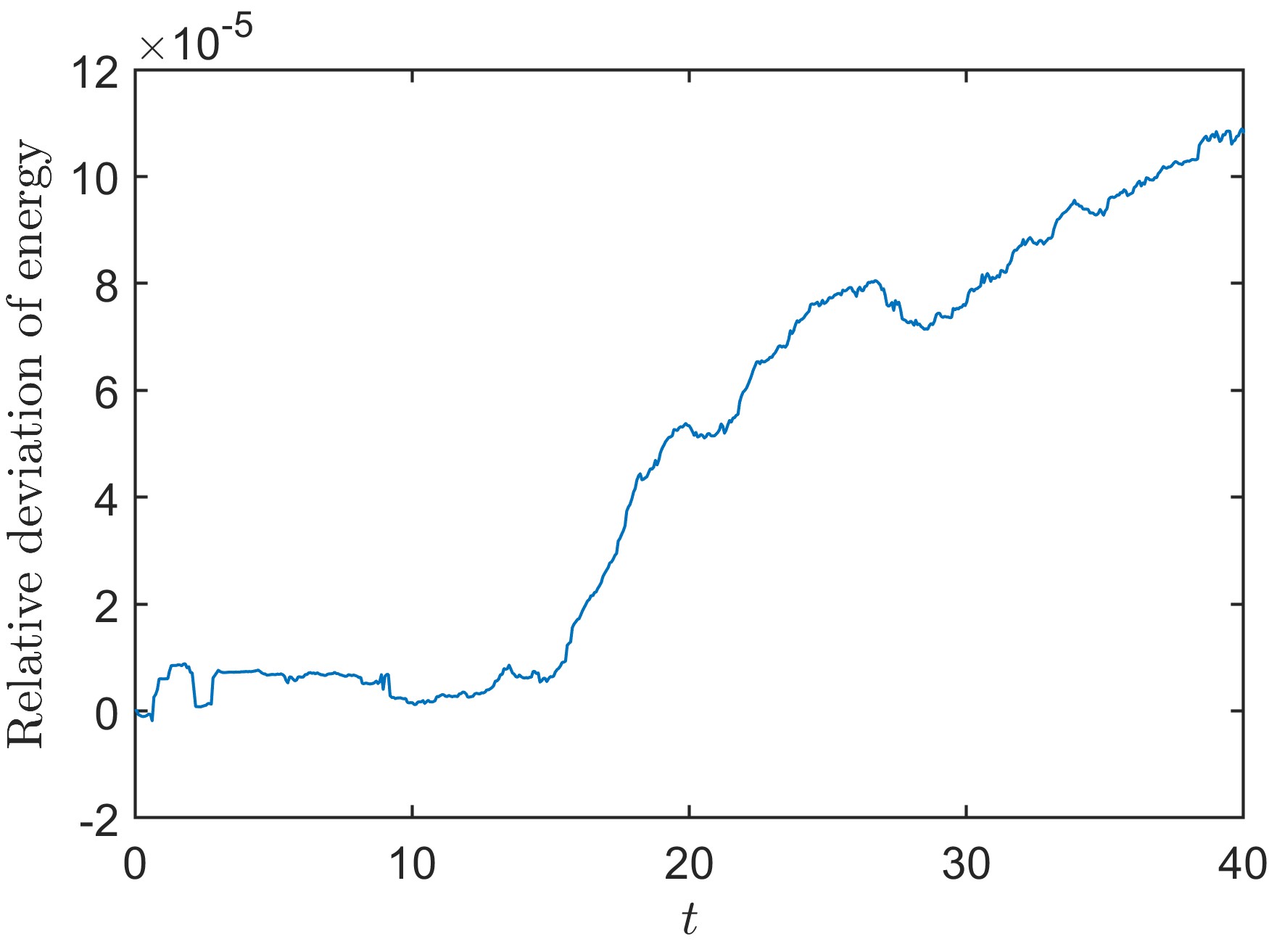}
  }

 \caption{(Bump-on-tail instability). Top left: time evolution of the $L^2$ norm of the electric field. Top right: performance of preserving mass. Bottom: performance of preserving momentum (left) and energy (right). 
 }
  \label{fig:mass_momentum_energy_BOT}
\end{figure}

\end{exa}

\section{Conclusion}\label{sec:conclusion}
In this paper, we combine the SL approach in time with cross approximation in space. Several novel ingredients are developed for this combination, including the local SL-FD solver, the new step-and-truncate strategy, the forward-characteristic tracing method for predicting information, and the new LoMaC method. The resulting SLAR methods are mass-conservative, high-order accurate in both time and space, and efficient in both computation and memory usage. An extensive set of numerical results are presented to demonstrate the effectiveness of the proposed SLAR methods. Future work includes generalizing the SLAR methods to high-dimensional  nonlinear multi-scale models with structure preserving properties.

\section{Acknowledgements} 

This work was partially supported Department of Energy DE-SC0023164 by the Multifaceted Mathematics Integrated Capability Centers (MMICCs) program
of DOE Office of Applied Scientific Computing Research (ASCR), and Air Force Office of Scientific Research (AFOSR) FA9550-24-1-0254 via the Multidisciplinary University Research Initiatives (MURI) Program. 
D. H., J.Q., N. Z. are also supported by NSF grant NSF-DMS-2111253, Air Force Office of Scientific Research FA9550-22-1-0390. 


\bibliographystyle{siam} 
\bibliography{Reference}

\end{document}